\documentclass[12pt]{article}
\usepackage{amsfonts}
\usepackage{mathrsfs}
\usepackage{amsthm}
\usepackage{graphics}
\usepackage{graphicx} 
\usepackage{amsmath} 
\usepackage{amssymb} 
\usepackage{enumerate}
\usepackage{epsfig}
\usepackage{multirow}
\usepackage{color}
\usepackage{cite}
\usepackage{tikz}

\usepackage{appendix}

\newtheorem{theorem}{\hskip\parindent Theorem}
\newtheorem{lemma}{\hskip\parindent Lemma}

\newtheorem{corollary}{\hskip\parindent Corollary}

\numberwithin{equation}{section}

\def\Ai{\mathrm {Ai}}
\def\Bi{\mathrm {Bi}}

\def\Im {\mathop{\rm Im}\nolimits}

\def\sign {\mathrm {sgn}}
\def\B {\mathrm{B}}

\begin{document}

\title{Real solutions of the first Painlev\'{e} equation with large initial data}

\author{Wen-Gao Long$^\dag$,\quad Yu-Tian Li$^\natural$,\quad Sai-Yu Liu$^\ddag$\quad and\quad Yu-Qiu Zhao$^\dag$\footnote{Corresponding author (Yu-Qiu Zhao).
 {\it{E-mail address:}} {stszyq@mail.sysu.edu.cn} } }

 \date{{\small{{\it{$^\dag$Department of Mathematics, Sun Yat-sen University, GuangZhou
510275, China}}\\
 {\it{$^\natural$Department of Mathematics, and Institute of Computational and Theoretical Studies,
Hong Kong Baptist University, Kowloon, Hong Kong}}\\
{\it{$^\ddag$School of Mathematics and Computation Science, Hunan University of Science and Technology, Xiangtan  411201, China}}
}}}
\maketitle
\begin{abstract}
We consider three special cases of the initial value problem of the first Painlev\'{e} equation (PI). Our approach is based on the method of {\it uniform asymptotics} introduced by Bassom, Clarkson, Law and McLeod. A rigorous proof of a property of the PI solutions on the negative real axis, recently revealed by Bender and Komijani, is given by approximating the Stokes multipliers. Moreover, we build more precise relation  between the large initial data of the PI solutions and their three different types of behavior  as the independent variable tends to negative infinity. In addition, some limiting form connection formulas are obtained.
\end{abstract}

\vspace{3mm}
{\it Keywords: }Connection formula; first Painlev\'{e} transcendent; uniform asymptotics; Airy function; modified Bessel function; Stokes multiplier

\vskip .2cm
{\it MSC2010:} 34M40, 33E17, 34A12, 34E05, 33C10

\section{Introduction and main results}
The first Painlev\'{e} equation
has the canonical form
\begin{equation}\label{PI equation}
\frac{d^{2}y}{dt^{2}}=6y^2+t.
\end{equation}
It is known  that there are two kinds of solutions of PI, behaving respectively as
$$y(t)\sim\sqrt{-\frac{t}{6}}\quad \text{and}\quad y(t)\sim-\sqrt{\frac{-t}{6}}$$
as $t\rightarrow-\infty$; see for example Bender and Orszag \cite{Bender-Orszag}. Refinements  have been obtained for the second case, such that  the solutions oscillate stably about the parabola $y=-\sqrt{\frac{-t}{6}}$, with
\begin{equation}\label{oscillate-behavior-t-infty}
y=-\left (-\frac{t}{6}\right )^{\frac{1}{2}}+d\,(-t)^{-\frac{1}{8}}\cos{\phi(t)}+\mathcal{O}\left (t^{-\frac{5}{8}}\right )~~\mbox{as}~~t\rightarrow-\infty,
\end{equation}
where
\begin{equation}\label{phi(t)}
\phi(t)=24^{\frac{1}{4}}\left[\frac{4}{5}(-t)^{\frac{5}{4}}-\frac{5}{8}d^2 \log(-t)+\theta\right],
\end{equation} and $d>0$ and $\theta$ are constants; cf. Kapaev \cite{AAKapaev-1988}, Joshi and Kruskal \cite{Joshi-Kruskal-1992},  and  Qin and Lu \cite{Qin-Lu-2008}, see also $\S 32.11$ of  the handbook \cite{NIST-handbook}.
Moreover, numerical analysis conducted in  Holmes and  Spence \cite{HS-1984}, Fornberg and Weideman \cite{Fornberg-Weideman} and Qin and Lu \cite{Qin-Lu-2008}  indicates that there exist  constants $\kappa_{1}<0$ and $\kappa_{2}>0$, such that all solutions of  \eqref{PI equation}  with $y(0)=0$ and $\kappa_{1}<y'(0)<\kappa_{2}$ belong to the second kind, while otherwise,  if $y'(0)>\kappa_{2}$ or $y'(0)<\kappa_{1}$, the solutions will blow up on the negative real axis.

Therefore, it is natural to consider the initial value problem of  \eqref{PI equation}  with   initial data
\begin{equation}\label{PI-equation-initial-problem}
y(0)=a,\quad y'(0)=b.
\end{equation}
The problem is how to determine the asymptotic behavior of $y(t)$ with a given real pair $(a, b)$. In particular, if $y(t)\sim-\sqrt{-t/6}$ as $t\rightarrow-\infty$, how to establish the relation  between the initial data (\ref{PI-equation-initial-problem}) and the parameters $d$ and $\theta$ in (\ref{oscillate-behavior-t-infty})-(\ref{phi(t)}). This is an open problem mentioned by  Clarkson in several occasions \cite{CPA2003, CPA2006}.

In fact, before Clarkson's open problem is being proposed,   investigations \cite{AAKapaev-1988, Kitaev-1995, HS-1984, Joshi-Kruskal-1992} have been made on the PI functions on the negative real axis. In \cite{HS-1984}, a boundary value problem for PI is studied by Holmes and  Spence, and it is shown that there are exactly three types of real solutions of PI equation. Later on, Kapaev \cite{AAKapaev-1988} obtained the asymptotic behavior of these solutions  as follows:
\begin{enumerate}
\item [(A)] a two-parameter family of solutions, oscillating about the parabola $y=-\sqrt{-t/6}$ and satisfying (\ref{oscillate-behavior-t-infty}) as $t\rightarrow-\infty$, where
    \begin{equation}\label{eq-parameter-d-theta}
  \left\{\begin{aligned}
  &24^{\frac{1}{4}}d^{2}=-\frac{1}{\pi}\log{|s_{0}|},\\
  &24^{\frac{1}{4}}\theta=-\arg{s_{3}}-24^{\frac{1}{4}}d^{2}\left(\frac{19}{8}\log{2}+\frac{5}{8}\log{3}\right)-\frac{\pi}{4}-\arg\Gamma\left(-i\frac{24^{\frac{1}{4}}}{2}d^{2}\right);
\end{aligned}\right.
\end{equation}
\item [(B)] a one-parameter family of solutions (termed {\it separatrix solutions}), satisfying
    \begin{equation}\label{eq-behavior-type-B}
      y(t)=\sqrt{\frac{-t}{6}}-\frac{h}{4\sqrt{\pi}}24^{-\frac{1}{8}}(-t)^{-\frac{1}{8}}\exp\left\{-\frac{4}{5}24^{\frac{1}{4}}(-t)^{\frac{5}{4}}\right\}+\mathcal{O}\left(|t|^{-\frac{5}{2}}\right)
    \end{equation}
    as $t\rightarrow-\infty$, where
    \begin{equation}\label{eq-parameter-h}
     h=s_{1}-s_{4};
    \end{equation}
\item [(C)] a two-parameter family of solutions, having infinitely many double poles on the negative real axis and satisfying
    \begin{equation}\label{eq-behavior-type-C}
      \frac{1}{y(t)+\sqrt{{-t}/{6}}}\sim \frac{\sqrt{6}}{2}\sin^{2}\left\{\frac{2}{5}24^{1/4}(-t)^{\frac{5}{4}}+\frac{5}{8}\rho\log(-t)+\sigma\right\}\quad \text{as}\quad t\rightarrow-\infty,
    \end{equation}
    where
    \begin{equation}\label{eq-parameter-rho-sigma}
    \left\{\begin{aligned}
    \rho&=\frac{1}{2\pi}\log(|s_{2}|^{2}-1)=\frac{1}{2\pi}\log(|1+s_{2}s_{3}|)=\frac{1}{2\pi}\log|s_{0}|,\\
    \sigma&=\frac{19}{8}\rho\log{2}+\frac{5}{8}\rho\log{3}+\frac{1}{2}\arg\Gamma\left(\frac{1}{2}-i\rho\right)-\frac{\pi}{4}+\frac{1}{2}\arg{s_{2}}.
    \end{aligned}\right.
    \end{equation}
\end{enumerate}
In the above formulas, $s_{k}$ for $k=0,1,2,3,4$  are the Stokes multipliers associated with the given solution; see the definition of Stokes multipliers in (\ref{eq-Stokes-matrices}) below. Moreover, using the isomonodromy condition, Kapaev \cite{AAKapaev-1988} also got the following necessary conditions for the Stokes multipliers of each solution type:
\begin{equation}\label{eq-condition-stokes-multipliers}
\begin{aligned}
1+s_{2}s_{3}>0 \quad \text{for type (A)},\\
1+s_{2}s_{3}=0 \quad \text{for type (B)},\\
1+s_{2}s_{3}<0 \quad \text{for type (C)}.
\end{aligned}
\end{equation}
 It should be noted that the solutions of types (A) and (B) may also have finite poles on the negative real axis; see \cite[Figures 1 and 2]{Bender-Komijani-2015}.

In the present  paper, we focus on three special cases of the initial value problem (\ref{PI-equation-initial-problem}) of the PI equation \eqref{PI equation}. More precisely, we shall establish the relation  between the asymptotic behavior  of the PI solution $y(t)$  as $t\to -\infty$ and the  real initial data $(y(0),y'(0))=(a,b)$,  in the following     cases:
\begin{equation*}
\begin{array}{ll}
   \text{(I)}& \text{fixed   $a$ and large positive (or, negative) $b$};\\
\text{(II)}& \text{fixed   $b$ and large negative $a$};\\
\text{(III)}& \text{fixed $b$ and large positive $a$}.
  \end{array}
\end{equation*}
The motivation of the first two cases comes from a recent investigation by Bender and Komijani \cite{Bender-Komijani-2015}, in which the unstable separatrix solutions of PI on the negative real axis are studied  numerically and analytically. For case (I), Bender and Komijani conclude that for any fixed initial value $y(0)$, there exists a sequence of initial slopes $y'(0)=b_{n}$ that give rise to separatrix solutions.   Figures 1 and 2 in \cite{Bender-Komijani-2015} indicate  an interesting phenomenon of the PI solutions as the initial slope varies. It seems that, when $b_{2n}<y'(0)<b_{2n+1}$, the solution passes through $n$ double poles and then oscillates stably about $-\sqrt{-t/6}$, while for $b_{2n-1}<y'(0)<b_{2n}$, the solution  has infinitely many  double poles. Moreover, they establish the asymptotic behavior
\begin{equation}\label{asymptotic-bn-En}
b_{n}\sim B n^{\frac{3}{5}}=2\left[\frac{\sqrt{3\pi}\Gamma(\frac{11}{6})n}{\Gamma(\frac{1}{3})}\right]^{\frac{3}{5}} \quad\text{as}\quad n\rightarrow +\infty
\end{equation}
by studying   the eigenvalue problem associated with the $\mathcal{PT}$-symmetric Hamiltonian $\hat{H}=\frac{1}{2}\hat{p}^{2}+2i \hat{x}^{3}$; see \cite{Bender-Komijani-2015}.
They have also obtained similar results for case (II).
However, part  of their argument  is based on numerical evidences (see \cite[p.8]{Bender-Komijani-2015}). To give a rigorous proof of their results and to build a precise relation  between the initial data and the large negative $t$ asymptotics, further analysis is  needed. Other than the above two cases, case (III) is essentially different, and has not been addressed  before, as far as we are aware of.



In the above   cases (I)-(III),
  we consider the initial value problem by  using the method of {\it uniform asymptotics}    introduced by Bassom, Clarkson, Law and McLeod \cite{APC}, and further developed by Wong and Zhang \cite{Wong-Zhang-2009-PIII, Wong-Zhang-2009-PIV}, Zeng and Zhao \cite{Zeng-Zhao-2015}, and Long, Zeng and Zhou \cite{Long-Zeng-Zhou-2017}. Our main results are stated as follows. Parts  of Theorems \ref{our theorem} and \ref{our theorem 2} are given  in \cite{Bender-Komijani-2015} and the rest are new.

First, for case (I), we shall prove the following results.
\begin{theorem}\label{our theorem}
Suppose that $y(t; a, b)$ is the real solution of (\ref{PI equation}) with initial data $(y(0),y'(0)) = (a,b)$, then for any fixed $a$, there exists an ascending sequence, $0<b_{1}<b_{2}<\cdots <b_{n}<\cdots$, such that $y(t; a,b_{n})$ are  separatrix solutions,  namely, of type (B), with $b_n$ possessing
  the large-$n$ asymptotic behavior (\ref{asymptotic-bn-En}), and the quantity $h$ in \eqref{eq-parameter-h}  having  the approximation
\begin{equation}\label{eq-limit-connection-h}
h:=h_n=(-1)^{n-1}2\exp\left\{\left(\frac{b_{n}}{2}\right)^{5/3}P_{0}\right\}\left(1+o(1)\right)  \quad \text{as}\quad n\rightarrow+\infty.
\end{equation}
Moreover,
\begin{enumerate}
\item [(i)] if $b_{2m}<b<b_{2m+1}$, $m=1,2,\cdots$, then $y(t; a,b)$ belongs to type (A) and
\begin{equation}\label{eq-limit-connection-d-large-b}
  \begin{aligned}
  24^{\frac{1}{4}}d^{2}=&\frac{1}{\pi}\left\{\left(\frac{b}{2}\right)^{\frac{5}{3}}P_{0}-\log\left[2\cos{\left(\left(\frac{b}{2}\right)^{\frac{5}{3}}Q_{0}- Q_{1} \right)}\right]+o(1)\right\},\\
  24^{\frac{1}{4}}\theta=&2\left(\frac{b}{2}\right)^{\frac{5}{3}}Q_{0}-2Q_{1}-24^{\frac{1}{4}}d^{2}\left(\frac{19}{8}\log{2}+\frac{5}{8}\log{3}\right)\\
  &+\frac{\pi}{4}-\arg\Gamma\left(-i\frac{24^{\frac{1}{4}}}{2}d^{2}\right)+o(1)
\end{aligned}
\end{equation}
as $b\rightarrow+\infty$;
\item [(ii)] if $b_{2m-1}<b<b_{2m}$, $m=1,2,\cdots$, then $y(t; a,b)$ belongs to type (C) and
\begin{equation}\label{eq-limit-connection-rho-sigma-case-I}
\begin{aligned}
  \rho&=-\frac{1}{2\pi}\left\{\left(\frac{b}{2}\right)^{\frac{5}{3}}P_{0}-\log\left[-2\cos{\left(\left(\frac{b}{2}\right)^{\frac{5}{3}}Q_{0}- Q_{1} \right)}\right]+o(1)\right\},\\
  \sigma&=\frac{19}{8}\rho\log{2}+\frac{5}{8}\rho\log{3}+\frac{1}{2}\arg\Gamma\left(\frac{1}{2}-i\rho\right)-\frac{\pi}{2}+\left(\frac{b}{2}\right)^{\frac{5}{3}}Q_{0}-Q_{1}+o(1)
\end{aligned}
\end{equation}
as $b\rightarrow+\infty$.
\end{enumerate}
In  (\ref{eq-limit-connection-h}), (\ref{eq-limit-connection-d-large-b}) and (\ref{eq-limit-connection-rho-sigma-case-I}), the constants $P_{0}, Q_{0}$ and $Q_{1}$ are expressible in terms of  the beta function as
\begin{equation}\label{eq-constants-thm-1}
\begin{aligned}
  P_{0}=\frac{6}{5}\B\left(\frac{1}{2},\frac{1}{3}\right),\quad
  Q_{0}=\frac{2\sqrt{3}}{5}\B\left(\frac{1}{2},\frac{1}{3}\right),\quad Q_{1}=\frac{\pi}{3}.
\end{aligned}
\end{equation}  There also exists a descending sequence, $\cdots<\tilde{b}_{n}<\cdots<\tilde{b}_{2}<\tilde{b}_{1}<0$,
having the similar properties as those of $\{b_{n}\}$.
\end{theorem}

For case (II), the following holds.
\begin{theorem}\label{our theorem 2}
Suppose that $y(t; a, b)$ is the real solution of (\ref{PI equation}) with initial data $(y(0),y'(0)) = (a,b)$, then for any real fixed $b$, there exists a negative descending sequence $0>a_{1}>a_{2}>\cdots >a_{n}>\cdots$, such that $y(t; a_{n}, b)$ are   separatrix solutions, namely, of type (B), with $a_{n}$ possessing  the asymptotic behavior
$$a_{n}\sim -\left[\frac{\sqrt{3\pi}\Gamma(11/6)}{\Gamma(1/3)}\right]^{\frac{2}{5}}n^{\frac{2}{5}}~~\text{as}~~n\rightarrow +\infty,$$
and the quantity $h$ in \eqref{eq-parameter-h}  having  the approximation
\begin{equation}\label{eq-limit-connection-h-case-II}
h:=h_n=(-1)^{n-1}2\exp\left\{\left(-a_{n}\right)^{5/2}P_{0}\right\}\left(1+o(1)\right) \quad \text{as}\quad n\rightarrow+\infty.
\end{equation} Moreover,
\begin{enumerate}
\item [(i)]  if $a_{2m+1}<a<a_{2m}$, $m=1,2,\cdots$, then $y(t; a, b)$ belongs to type (A) and
\begin{equation}\label{limit-connection-formula-d-large-a}
\begin{aligned}
24^{\frac{1}{4}}d^{2}=&\frac{1}{\pi}\left\{(-a)^{\frac{5}{2}}P_{0}-\log\left[2\cos{\left((-a)^{\frac{5}{2}}Q_{0}\right)}\right]+o(1)\right\},\\
24^{\frac{1}{4}}\theta=&2\left(-a\right)^{\frac{5}{2}}Q_{0}-24^{\frac{1}{4}}d^{2}\left(\frac{19}{8}\log{2}+\frac{5}{8}\log{3}\right)\\
  &+\frac{\pi}{4}-\arg\Gamma\left(-i\frac{24^{\frac{1}{4}}}{2}d^{2}\right)+o(1)
\end{aligned}
\end{equation}
as $a\rightarrow-\infty$;
\item [(ii)] if $a_{2m}<a<a_{2m-1}$, $m=1,2,\cdots$, then $y(t; a, b)$ belongs to type (C) and
\begin{equation}\label{eq-limit-connection-rho-sigma-case-II}
\begin{aligned}
  \rho&=-\frac{1}{2\pi}\left\{\left(-a\right)^{\frac{5}{2}}P_{0}-\log\left[-2\cos{\left(\left(-a\right)^{\frac{5}{2}}Q_{0}\right)}\right]+o(1)\right\},\\
  \sigma&=\frac{19}{8}\rho\log{2}+\frac{5}{8}\rho\log{3}+\frac{1}{2}\arg\Gamma\left(\frac{1}{2}-i\rho\right)-\frac{\pi}{2}+\left(-a\right)^{\frac{5}{2}}Q_{0}+o(1)
\end{aligned}
\end{equation}
as $a\rightarrow-\infty$.
\end{enumerate}
In   (\ref{eq-limit-connection-h-case-II}), (\ref{limit-connection-formula-d-large-a}) and (\ref{eq-limit-connection-rho-sigma-case-II}), the constants $P_{0}$ and $Q_{0}$ are the same as the ones in (\ref{eq-constants-thm-1}).
\end{theorem}

Case (III) is in a sense different from the first two cases. We have   the following result.
\begin{theorem}\label{our-theorem-3}
  Suppose that $y(t; a, b)$ is the real solution of (\ref{PI equation}) with initial data (\ref{PI-equation-initial-problem}). For any fixed $b$, there exists a positive number $M$ (depending on $b$) such that, if $a>M$, then $y(t; a, b)$ belongs to type (C), and
\begin{equation}\label{eq-limit-connection-rho-sigma-case-III}
\begin{aligned}
  \rho&=\frac{1}{\pi}\left [a^{\frac{5}{2}}H_{0}+o(1)\right ],\\
  \sigma&=\frac{19}{8}\rho\log{2}+\frac{5}{8}\rho\log{3}+\frac{1}{2}\arg\Gamma\left(\frac{1}{2}-i\rho\right)-\frac{\sqrt{3}}{2}a^{\frac{5}{2}}H_{0}
  +o(1),
\end{aligned}
\end{equation}
where $H_{0}=\frac{2}{5}\B\left(\frac{1}{2},\frac{1}{6}\right)$.
\end{theorem}\vskip 1cm

For the sake of  convenience, we   recall some important concepts in the isomonodromy theory for the first Painlev\'{e} transcendents. First, one of the Lax pairs for the PI equation is given as follows (see \cite{Kapaev-Kitaev-1993}):
\begin{equation}\label{lax pair-I}
\left\{\begin{aligned}
\frac{\partial\Psi}{\partial\lambda}&=\left\{(4\lambda^4+t+2y^2)\sigma_{3}-i(4y\lambda^2+t+2y^2)\sigma_{2}-(2y_{t}\lambda+\frac{1}{2\lambda})\sigma_{1}\right\}\Psi=\mathcal{A}(\lambda)\Psi,\\
\frac{\partial\Psi}{\partial t}&=\left\{(\lambda+\frac{y}{\lambda})\sigma_{3}-\frac{iy}{\lambda}\sigma_{2}\right\}\Psi=\mathcal{B}(\lambda)\Psi,
\end{aligned}\right.
\end{equation}
where
$$\sigma_{1}=\left(\begin{matrix}0&1\\1&0\end{matrix}\right),\quad \sigma_{2}=\left(\begin{matrix}0&-i\\i&0\end{matrix}\right),\quad \sigma_{3}=\left(\begin{matrix}1&0\\0&-1\end{matrix}\right)$$
are the Pauli matrices and $y_{t}=\frac{dy}{dt}$. Compatibility of the above system implies that $y=y(t)$ satisfies the first Painlev\'{e} equation (\ref{PI equation}). Under the transformation
\begin{equation}\label{eq-transform-canonical solution}
\Phi(\lambda)=\lambda^{\frac{1}{4}\sigma_{3}}\frac{\sigma_{3}+\sigma_{1}}{\sqrt{2}}\Psi(\sqrt{\lambda}),
\end{equation}
the first equation of (\ref{lax pair-I}) becomes
\begin{equation}\label{eq-fold-Lax-pair}
\frac{\partial\Phi}{\partial\lambda}=\left(\begin{matrix}y_{t}&2\lambda^{2}+2y\lambda-t+2y^2\\2(\lambda-y)&-y_{t}\end{matrix}\right)\Phi.
\end{equation}
Following \cite{Kapaev-Kitaev-1993} (see also \cite{AAKapaev-2004}), the only singularity of equation (\ref{eq-fold-Lax-pair}) is the  irregular singular point at $\lambda=\infty$, and there exist canonical solutions $\Phi_{k}(\lambda)$, $k\in\mathbb{Z}$, of (\ref{eq-fold-Lax-pair}) with the following asymptotic expansion
\begin{equation}\label{eq-canonical-solutions}
\Phi_{k}(\lambda,t)=\lambda^{\frac{1}{4}\sigma_{3}}\frac{\sigma_{3}+\sigma_{1}}{\sqrt{2}}\left(I+\frac{\mathcal{H}}{\lambda}+\mathcal{O}\left(\frac{1}{\lambda^2}\right)\right)e^{(\frac{4}{5}\lambda^{\frac{5}{2}}+t\lambda^{\frac{1}{2}})\sigma_{3}},~~\lambda\rightarrow\infty,~~\lambda\in\Omega_{k},
\end{equation}
where $\mathcal{H}=-(\frac{1}{2}y_{t}^2-2y^3-ty)\sigma_{3}$, and the canonical sectors are
$$\Omega_{k}=\left\{\lambda\in\mathbb{C}:~\arg \lambda\in \left(-\frac{3\pi}{5}+\frac{2k\pi}{5},\frac{\pi}{5}+\frac{2k\pi}{5}\right)\right\}, \qquad k\in\mathbb{Z}.$$
These canonical solutions are related by
\begin{equation}\label{eq-Stokes-matrices}
\Phi_{k+1}=\Phi_{k}S_{k},\quad S_{2k-1}=\left(\begin{matrix}1&s_{2k-1}\\0&1\end{matrix}\right),\quad S_{2k}=\left(\begin{matrix}1&0\\s_{2k}&1\end{matrix}\right),
\end{equation}
where $s_{k}$   are called Stokes multipliers, and   independent of $\lambda$ and $t$ according to the isomonodromy condition. The Stokes multipliers are   subject to the constraints
\begin{equation}\label{eq-constraints-stokes-multipliers}
s_{k+5}=s_{k}\quad \text{and}\quad s_{k}=i(1+s_{k+2}s_{k+3}),~~k\in\mathbb{Z}.
\end{equation}
Moreover,  regarding $s_{k}$ as functions of $(t,y(t),y'(t))$,   they also satisfy \cite[p.1687, (13)]{AAKapaev-1988}
\begin{equation}\label{eq-sk-s-k-relation}
s_{k}\left (t,y(t),y'(t)\right)=-\overline{s_{-k}\left (\bar{t},\overline{y(t)},\overline{y'(t)}\right )}, \quad k\in\mathbb{Z},
\end{equation}
where  $\bar\zeta$ stands for the complex    conjugate of a complex number   $\zeta$.  From \eqref{eq-constraints-stokes-multipliers}, it is readily seen that, in general, two of the Stokes multipliers determine all others. The derivation  of (\ref{eq-canonical-solutions}), (\ref{eq-Stokes-matrices}) and (\ref{eq-constraints-stokes-multipliers}), and more details about the Lax pairs, are referred to in \cite{FAS-2006}.

The rest of the paper is arranged as follows. In Sec.\;\ref{sec:proof-main-result}, we prove the main theorems, assuming the validity of Lemmas \ref{Thm-s0-large-slope}, \ref{Thm-s0-large-value} and \ref{Thm-large-initial-value-postive}. Next, in Sec.\,\ref{sec:proof-lemma},   we apply the method of {\it uniform asymptotics} to calculate the Stokes multipliers when $t=0$,  and hence to   prove the three lemmas.  In the last section, Sec.\,\ref{sec:discussion}, a discussion is  provided  on  prospective studies and possible difficulties  in the general case of Clarkson's open problem.  Some technical details are put in Appendices A, B, C and D to clarify the derivation.

\section{Proof of the main results}\label{sec:proof-main-result}
Because the monodromy data $\{s_0, s_1, s_2, s_3, s_4\}$ and the solutions of PI have a one-to-one correspondence,
a general idea to solve connection problems is to calculate the Stokes multipliers of a specific solution  in the two specific situations to be connected. In the initial data case,  this  means   to calculate all $s_{k}$  as   $t\rightarrow-\infty$ and  at $t=0$. When $t\rightarrow-\infty$, as stated above in (\ref{eq-parameter-d-theta}), (\ref{eq-parameter-h}) and (\ref{eq-parameter-rho-sigma}), the Stokes multipliers have been derived by Kapaev \cite{AAKapaev-1988}. 
When $t=0$,   it is much difficult to find the exact values of $s_{k}$. However, inspired  by the ideas in Sibuya \cite{Sibuya-1967},  we are able to obtain  their  asymptotic approximations  in the special cases considered here, as a step forward. These approximations, together with (\ref{eq-parameter-d-theta}), (\ref{eq-parameter-h}), (\ref{eq-parameter-rho-sigma}) and (\ref{eq-condition-stokes-multipliers}), suffice to prove our theorems.

First, the following lemma is crucial to prove Theorem \ref{our theorem}.
\begin{lemma}
\label{Thm-s0-large-slope}
For any fixed $a$ and large positive (or negative) $b$, the asymptotic behaviors of the Stokes multipliers, corresponding to $y(t; a, b)$, are given by
\begin{equation}\label{s0-large slope}
\begin{aligned}
s_{0}&=2i\exp\left\{-\left(\frac{|b|}{2}\right)^{\frac{5}{3}} (E_{0}+F_{0})\right\}\left\{\cos{\left[\left(\frac{|b|}{2}\right)^{\frac{5}{3}}Q_{0}-\sign(b) \cdot Q_{1}\right]+o(1)}\right\},\\
s_{1}&=i\exp\left\{2\left(\frac{|b|}{2}\right)^{\frac{5}{3}}E_{0}-2\cdot\sign(b)\cdot E_{1}+o(1)\right\}=-\overline{s_{4}},\\
s_{3}&=-i\exp\left\{2\left(\frac{|b|}{2}\right)^{\frac{5}{3}}(E_{0}-F_{0})-2\cdot \sign(b)\cdot(E_{1}-F_{1})+o(1)\right\}=-\overline{s_{2}}
\end{aligned}
\end{equation}
as $|b|\rightarrow+\infty$, where $\sign(b)$ is the sign of $b$ and
\begin{equation}\label{eq-constants-lemma-1}
\begin{aligned}
E_{0}&=\overline{F_{0}}=\frac{3}{5}\B\left(\frac{1}{2},\frac{1}{3}\right)-\frac{\sqrt{3}i}{5}\B\left(\frac{1}{2},\frac{1}{3}\right),\quad E_{1}=-F_{1}=-\frac{\pi i}{6},\\[0.2cm]
Q_{0}&=i(E_{0}-F_{0})=\frac{2\sqrt{3}}{5}\B\left(\frac{1}{2},\frac{1}{3}\right),\quad Q_{1}=i(E_{1}-F_{1})=\frac{\pi}{3}.
\end{aligned}
\end{equation}
\vskip .4 cm

\end{lemma}
An immediate consequence of (\ref{s0-large slope}) is stated as follows.
\begin{corollary}
For any fixed $a$, there exists a positive sequence $\{b_{n}\}$ and a negative sequence $\{\tilde{b}_{n}\}$, $n=1,2,\cdots$,  with
 \begin{equation}\label{asymptotic-bn}
 b_{n}\sim 2\left(\frac{n\pi-\pi/2+Q_{1}}{Q_{0}}\right)^{\frac{3}{5}}\quad \text{and} \qquad \tilde{b}_{n}\sim -2\left(\frac{n\pi-\pi/2-Q_{1}}{Q_{0}}\right)^{\frac{3}{5}},
 \end{equation}
as $n\rightarrow\infty$, such that the Stokes multipliers, corresponding to $y(t; a,b_{n})$ (or $y(t; a, \tilde{b}_{n})$), satisfy
$$s_{0}=0,\quad s_{3}=s_{2}=s_{1}+s_{4}=i.$$
\end{corollary}\vskip .4 cm

\begin{table}[h]
  \centering
\setlength{\belowcaptionskip}{-0.3cm}
  \begin{tabular}{|c|c|c|c|c|c|c|c|c|}
  \hline
  $n$&1&2&3&4&5&&10&11\\
  \hline
  the true values of $b_{n}$&1.8518&3.0040&3.9052&4.6834&5.3831&$\cdots$&8.2449&8.7383\\
  \hline
  $B n^{\frac{3}{5}}$ &2.0922&3.1711&4.0445&4.8065&5.4951&$\cdots$&8.3290&8.8192\\
  \hline
  $2\left(\frac{n\pi-\pi/2+Q_{1}}{Q_{0}}\right)^{3/5}$&1.8754&3.0098&3.9081&4.6853&5.3844&$\cdots$ &8.2454 & 8.7388\\
  \hline
  \end{tabular}
  \caption{\small When $a=0$, the comparison of the approximate values (\ref{asymptotic-bn-En}), (\ref{asymptotic-bn}) and the true values of $b_{n}$, using five significant digits. The true values of $b_{n}$ are taken from \cite{Bender-Komijani-2015}.}\label{Table-bn asymptotic}
\end{table}

\textbf{Proof of Theorem \ref{our theorem}.} We only give  the proof when $b>0$. When $b<0$, the argument is the same except for some minor justification of signs; see (\ref{s0-large slope}). According to (\ref{eq-condition-stokes-multipliers}), see also \cite[Theorems 2.1 and 2.2]{AAKapaev-2004},   regarding $s_{0}$ as a function of $a$ and $b$, then $s_{0}(a, b_{n})=0$ implies that the solutions $y(t; a, b_{n})$ all belong to type (B). Moreover, (\ref{asymptotic-bn}) implies (\ref{asymptotic-bn-En}), and is an improved approximation of $b_{n}$. Although both  formulas are only applicable   for large $n$, (\ref{asymptotic-bn}) has a better accuracy for small $n$ when $a=0$ as compared with (\ref{asymptotic-bn-En});  see Table \ref{Table-bn asymptotic}. For fixed $a$, from (\ref{s0-large slope}), it is obvious that the sequence $\{b_{n}\}$ can be chosen such that
\begin{equation*}
1+s_{2}s_{3}=\Im s_{0}(a,b)
\left\{\begin{aligned}
&>0, \quad \text{if}\quad b_{2m}<b<b_{2m+1},\\
&<0, \quad \text{if}\quad b_{2m-1}<b<b_{2m},
\end{aligned}\right.\quad m=1,2,\cdots.
\end{equation*}
It means that if $b$ varies in the above two sequences of  open  intervals, the PI solutions $y(t; a, b)$ correspond alternatively to type (A) and type (C) solutions.

To get (\ref{eq-limit-connection-h}), we   calculate the leading asymptotic behavior of $s_{1}-s_{4}$ as $b\rightarrow+\infty$. In fact, as a consequence of (\ref{s0-large slope}), we get
\begin{equation}\label{eq-s1-s4-behavior}
\begin{aligned}
  s_{1}-s_{4}&=i\exp\left\{2\left(\frac{b}{2}\right)^{\frac{5}{3}} E_{0}-2 E_{1}+o(1)\right\}-i\exp\left\{2\left(\frac{b}{2}\right)^{\frac{5}{3}}F_{0}-2F_{1}+o(1)\right\}\\
             &=2\exp\left\{\left(\frac{b}{2}\right)^{\frac{5}{3}}(E_{0}+F_{0})\right\}\left\{\sin\left(i\left(\frac{b}{2}\right)^{\frac{5}{3}}(E_{0}-F_{0})-i (E_{1}-F_{1})\right)+o(1)\right\}\\
             &=2\exp\left\{\left(\frac{b}{2}\right)^{\frac{5}{3}}P_{0}\right\}\left\{\sin\left(\left(\frac{b}{2}\right)^{\frac{5}{3}} Q_{0}-Q_{1}\right)+o(1)\right\}
\end{aligned}
\end{equation}
as $b\rightarrow+\infty$. Here we have set $P_{0}=E_{0}+F_{0}$. Putting $b=b_{n}$, and noting that $\left(\frac{b_{n}}{2}\right)^{5/3} Q_{0}-Q_{1}=n\pi-\frac{\pi}{2}+o(1)$ as $n\rightarrow\infty$, we obtain (\ref{eq-limit-connection-h}).

Next, to obtain the limiting form connection formulas in (\ref{eq-limit-connection-d-large-b}), we only need to calculate $|s_{0}|$ and $\arg s_{3}$. According to (\ref{s0-large slope}) and (\ref{eq-constants-lemma-1}), one has
\begin{equation}\label{eq-s0-mo}
\begin{aligned}
  &|s_{0}|=2\exp\left\{-\left(\frac{b}{2}\right)^{5/3}(E_{0}+F_{0})\right\}\left\{\cos{\left(\left(\frac{b}{2}\right)^{5/3}Q_{0}- Q_{1}\right)}+o(1)\right\},\\
  &\arg{s_{3}}=-\frac{\pi}{2}-\left(\frac{b}{2}\right)^{5/3}Q_{0}+Q_{1}+o(1)
\end{aligned}
\end{equation}
as $b \rightarrow+\infty$. Substituting (\ref{eq-s0-mo}) into (\ref{eq-parameter-d-theta}) and noting that $P_{0}=E_{0}+F_{0}$, one may immediately get (\ref{eq-limit-connection-d-large-b}).
Finally, to obtain (\ref{eq-limit-connection-rho-sigma-case-I}), we   calculate $|s_{0}|$ and $\arg s_{2}$,  which can also be derived from (\ref{s0-large slope}). This completes the proof of Theorem \ref{our theorem}.\qed
\\

Similarly, based on  the following result, we can prove Theorem \ref{our theorem 2}.
\begin{lemma}\label{Thm-s0-large-value}
For fixed $b$ and large negative $a$, the asymptotic behaviors of the Stokes multipliers, corresponding to $y(t; a, b)$, are given by
\begin{equation}
\begin{aligned}
s_{0}&=2i\exp\left\{-(-a)^{5/2}(E_{0}+F_{0})\right\}\left\{\cos{\left[(-a)^{5/2}Q_{0}\right]+o(1)}\right\},\\
s_{1}&=i\exp\left\{2(-a)^{5/2}E_{0}+o(1)\right\}=-\overline{s_{4}},\\ s_{3}&=i\exp\left\{2(-a)^{5/2}(E_{0}-F_{0})+o(1)\right\}=-\overline{s_{2}}
\end{aligned}
\end{equation}
as $a\rightarrow-\infty$, where $Q_{0}, E_{0}$ and $F_{0}$ are the same constants as the ones in \eqref{eq-constants-lemma-1}.
\end{lemma}

\begin{corollary}
There exists a negative discrete set $\{a_{n}\}$, $n=1,2,\cdots$, with
\begin{equation}
a_{n}\sim-\left(\frac{n\pi-\frac{\pi}{2}}{Q_{0}}\right)^{\frac{2}{5}}\quad \text{ as } n\rightarrow\infty,
\end{equation}
such that the Stokes multipliers corresponding to $y(t; a_{n}, b)$ are
$$s_{0}=0,\quad s_{3}=s_{2}=s_{1}+s_{4}=i.$$
\end{corollary}

Furthermore, Theorem \ref{our-theorem-3} is an immediate consequence of a combination of Lemma \ref{Thm-large-initial-value-postive} and (\ref{eq-condition-stokes-multipliers}).
\begin{lemma}\label{Thm-large-initial-value-postive}
For fixed $b$ and large positive $a$, the asymptotic behaviors of the Stokes multipliers, corresponding to $y(t; a, b)$, are given by
\begin{equation}\label{eq-s0-large-positive-a}
\begin{aligned}
  s_{0}&=-i\exp\left\{2 a^{\frac{5}{2}} H_{0}+o(1)\right\},\\
  s_{1}&=i\exp\left\{2 a^{\frac{5}{2}} G_{0}+o(1)\right\}=-\overline{s_{4}},\\
  s_{3}&=i\exp\left\{2 a^{\frac{5}{2}}(G_{0}+H_{0})+o(1)\right\}=-\overline{s_{2}}
\end{aligned}
\end{equation}
as $a\rightarrow+\infty$, where $H_{0}=\frac{2}{5}\B\left(\frac{1}{6},\frac{1}{2}\right)$ and $G_{0}=-\frac{1}{5}\B\left(\frac{1}{2},\frac{1}{6}\right)+\frac{\sqrt{3}i}{5}\B\left(\frac{1}{2},\frac{1}{6}\right)=e^{\frac{2\pi i}{3}}H_{0}$.
\end{lemma}

We leave  the proof of Lemmas \ref{Thm-s0-large-slope}, \ref{Thm-s0-large-value} and \ref{Thm-large-initial-value-postive} to the next section.

\section{Uniform asymptotics and the proofs of the  lemmas}\label{sec:proof-lemma}
In this section, we are going to prove Lemmas \ref{Thm-s0-large-slope}, \ref{Thm-s0-large-value} and \ref{Thm-large-initial-value-postive}, using the method of {\it uniform asymptotics} \cite{APC}. The method consists of two main steps. The first step is to transform the Lax pair equation (\ref{eq-fold-Lax-pair}) into a second-order Schr\"{o}dinger equation, and to approximate the solutions of this equation with well known special functions. Indeed,   denoting
$$\Phi(\lambda)=\left(\begin{matrix}\phi_{1}\\\phi_{2}\end{matrix}\right),$$
and defining $Y(\lambda,t)=(2(\lambda-y))^{-\frac{1}{2}}\phi_{2}$, we have
\begin{equation}\label{Schrodinger-equation-t-general}
\frac{d^{2}Y}{d\lambda^{2}}=\left[y_{t}^{2}+4\lambda^{3}+2\lambda t-2y t-4y^{3}-\frac{y_{t}}{\lambda-y}+\frac{3}{4}\frac{1}{(\lambda-y)^2}\right]Y(\lambda,t).
\end{equation}
When $t=0$, equation (\ref{Schrodinger-equation-t-general}) is simplified to
\begin{equation}\label{Schrodinger-equation}
\frac{d^{2}Y}{d\lambda^{2}}=\left[4\lambda^{3}+b^{2}-4a^{3}-\frac{b}{\lambda-a}+\frac{3}{4}\frac{1}{(\lambda-a)^2}\right]Y(\lambda).
\end{equation}
One can regard (\ref{Schrodinger-equation}) as either a scalar or a $1\times 2$ vector equation. We will see that in all   cases (I), (II) and (III) considered in this paper, the solutions of this Schr\"{o}dinger equation can be approximated by certain special functions. In the second step, we can therefore use the  known Stokes phenomena of these special functions. In each case,  we shall calculate the Stokes multipliers of $Y$, and then calculate those of $\Phi$. We carry out a case by case analysis to complete the steps.

\subsection{Case (I): fixed $a$ and large positive (negative) $b$}
We may assume that  $b>0$. The argument for the case when $b<0$ is the same, and hence omitted here.

With the scaling $\lambda=\xi^{\frac{2}{5}}\eta$, $b=2\xi^{\frac{3}{5}}$ as $\xi\rightarrow+\infty$, equation (\ref{Schrodinger-equation}) is reduced to
\begin{equation}\label{eq-shrodinger-y10-large}
\begin{split}
\frac{d^{2}Y}{d\eta^{2}}&=\xi^{2}\left[4(\eta+1)(\eta-e^{\frac{\pi i}{3}})(\eta-e^{\frac{-\pi i}{3}})-\frac{4a^2}{\xi^{\frac{6}{5}}}-\frac{2\xi^{-\frac{3}{5}}}{\xi^{\frac{2}{5}}\eta-a}+\frac{3\xi^{-\frac{6}{5}}}{4(\xi^{\frac{2}{5}}\eta-a)^{2}}\right]Y\\
&:=\xi^{2}F(\eta,\xi)Y,
\end{split}
\end{equation}
where
\begin{equation}\label{eq-F(eta,xi)}
  F(\eta,\xi)=4(\eta+1)(\eta-e^{\frac{\pi i}{3}})(\eta-e^{\frac{-\pi i}{3}})-\frac{2}{\xi\eta}+g(\eta,\xi),
\end{equation}
and $g(\eta,\xi)=\mathcal{O}\left(\xi^{-\frac{6}{5}}\right)$ as $\xi\rightarrow+\infty$, uniformly for all $\eta$ bounded away from $\eta=0$. Obviously, there are three simple turning points, say $\eta_{j}$, $j=0,1,2$. Those  are the  zeros of $F(\eta,\xi)$ near $-1, e^{\frac{\pi i}{3}}$ and $e^{-\frac{\pi i}{3}}$ respectively. For convenience, we denote $\alpha=e^{\frac{\pi i}{3}}$ and $\beta=e^{-\frac{\pi i}{3}}$. A straightforward  calculation shows that $\eta_{1}-\alpha\sim-\frac{1}{6\xi}$ and $\eta_{2}-\beta\sim-\frac{1}{6\xi}$  as $\xi\rightarrow +\infty$.

\begin{figure}[h]
  \centering
  \includegraphics[width=6cm, scale=1]{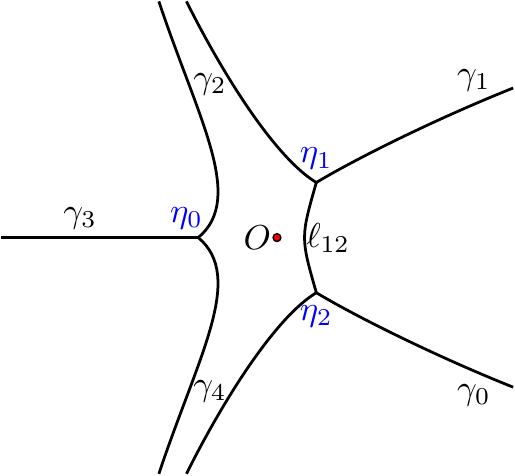}\\
  \caption{Stokes geometry of $F(\eta,\xi)d\eta^{2}$.}\label{Figure-stokes complex-large slope}
\end{figure}

According to \cite{EGS-2008}, the limiting state of the Stokes geometry of the quadratic form $F(\eta,\xi)d\eta^2$ as $\xi\rightarrow+\infty$ is described in Figure \ref{Figure-stokes complex-large slope}.
Therefore, following the main ideas in \cite{APC}, we can approximate the solutions of (\ref{eq-shrodinger-y10-large}) via the Airy functions. To do so, we  define two conformal mappings $\zeta(\eta)$ and $\omega(\eta)$ by
\begin{equation}\label{zeta-define}
\int_{0}^{\zeta}s^{\frac{1}{2}}ds=\int_{\eta_{1}}^{\eta}F(s,\xi)^{\frac{1}{2}}ds
\end{equation}
and
\begin{equation}\label{omega-define}
\int_{0}^{\omega}s^{\frac{1}{2}}ds=\int_{\eta_{2}}^{\eta}F(s,\xi)^{\frac{1}{2}}ds,
\end{equation}
respectively from neighborhoods of $\eta=\eta_1$ and $\eta=\eta_2$ to the origin. In the present paper, we take the principal branches for all the square roots. Then the conformality can be extended to the Stokes curves, and  the following lemma is a consequence of  \cite[Theorem 2]{APC}.

\begin{lemma}\label{lem-uniform-airy-large-slope}
There are constants $C_{1}, C_{2}$ and $\tilde{C}_{1}, \tilde{C}_{2}$, depending on $\xi$, such that
\begin{equation}
Y=\left(\frac{\zeta}{F(\eta,\xi)}\right)^{\frac{1}{4}}\left\{C_{1}\left[1+o(1)\right]\Ai(\xi^{\frac{2}{3}}\zeta)+C_{2}\left[1+o(1)\right]\Bi(\xi^{\frac{2}{3}}\zeta)\right\},\quad \xi\rightarrow+\infty
\end{equation}
uniformly for $\eta$ on any two adjacent Stokes lines emanating from $\eta_{1}$; and
\begin{equation}
Y=\left(\frac{\omega}{F(\eta,\xi)}\right)^{\frac{1}{4}}\left\{\tilde{C}_{1}\left[1+o(1)\right]\Ai(\xi^{\frac{2}{3}}\omega)+\tilde{C}_{2}\left[1+o(1)\right]\Bi(\xi^{\frac{2}{3}}\omega)\right\},\quad \xi\rightarrow+\infty
\end{equation}
uniformly for $\eta$ on any two adjacent Stokes lines emanating from $\eta_{2}$.
\end{lemma}
Moreover, we get the asymptotic behavior of $\zeta(\eta)$ and $\omega(\eta)$ as $\xi, |\eta|\rightarrow+\infty$.
\begin{lemma}\label{lemma-zeta-omega-eta-infty-relation}
\begin{equation}\label{zeta-eta-infty-relation}
\frac{2}{3}\zeta^{\frac{3}{2}}=\frac{4}{5}\eta^{\frac{5}{2}}+E_{0}-\frac{E_{1}}{\xi}+o\left(\frac{1}{\xi}\right),
\qquad \arg\eta\in \left(-\frac{\pi}{5},\pi\right),
\end{equation}
and
\begin{equation}\label{omega-eta-infty-relation}
\frac{2}{3}\omega^{\frac{3}{2}}=\frac{4}{5}\eta^{\frac{5}{2}}+F_{0}-\frac{F_{1}}{\xi}+o\left(\frac{1}{\xi}\right),
\qquad \arg\eta\in \left(-\pi,\frac{\pi}{5}\right)
\end{equation}
as $\xi,|\eta|\rightarrow +\infty$ with $|\eta|\gg\xi^{2}$, where $E_{1}, F_{1}, E_{0}$ and $F_{0}$ are constants stated in (\ref{eq-constants-lemma-1}).
\end{lemma}
The proof of Lemma \ref{lemma-zeta-omega-eta-infty-relation} is left to Appendix A. Now we turn to the proof of Lemma \ref{Thm-s0-large-slope}, assuming the validity of Lemma \ref{lemma-zeta-omega-eta-infty-relation}.\vskip .4cm

\textbf{Proof of Lemma \ref{Thm-s0-large-slope}.}
According to \cite{APC}, in order to calculate $s_{0}$, we should know the uniform asymptotic behaviors of $Y$ on the two Stokes lines   tending  to infinity with $\arg\eta\sim\pm\frac{\pi}{5}$. However, one may find that, in Figure \ref{Figure-stokes complex-large slope}, the two adjacent Stokes lines with $\arg\eta\sim \pm \frac{\pi}{5}$, emanate from two different turning points. Hence we should build the relations  between $(\Ai(\xi^{\frac{2}{3}}\zeta),\Bi(\xi^{\frac{2}{3}}\zeta))$ and $(\Ai(\xi^{\frac{2}{3}}\omega),\Bi(\xi^{\frac{2}{3}}\omega))$. Actually, it can be done by matching them on the Stokes line $\ell_{12}$ joining $\eta_{1}$ and $\eta_{2}$. According to \cite[(9.2.12) and (9.7.5)]{NIST-handbook}, one has
\begin{equation}\label{eq-behavior-Ai}
\left\{\begin{aligned}
\Ai(z)\sim&\frac{1}{2\sqrt{\pi}}z^{-\frac{1}{4}}e^{-\frac{2}{3}z^{\frac{3}{2}}}, &\arg z\in \left(-\pi,\pi\right);\\
\Ai(z)\sim&\frac{1}{2\sqrt{\pi}}z^{-\frac{1}{4}}e^{-\frac{2}{3}z^{\frac{3}{2}}}+\frac{i}{2\sqrt{\pi}}z^{-\frac{1}{4}}e^{\frac{2}{3}z^{\frac{3}{2}}}, &\arg z\in\left(\frac{\pi}{3},\frac{5\pi}{3}\right);\\
\Ai(z)\sim&\frac{1}{2\sqrt{\pi}}z^{-\frac{1}{4}}e^{-\frac{2}{3}z^{\frac{3}{2}}}-\frac{i}{2\sqrt{\pi}}z^{-\frac{1}{4}}e^{\frac{2}{3}z^{\frac{3}{2}}}, &\arg z\in\left(\frac{-5\pi}{3},\frac{-\pi}{3}\right)
\end{aligned}\right.
\end{equation}
as $|z|\rightarrow\infty$, and the uniform asymptotic behavior  of $\Bi(z)$ for $\arg z\sim\pm\frac{\pi}{3}$ can be evaluated by applying \cite[(9.2.10)]{NIST-handbook}, namely,
\begin{equation}\label{eq-Bi(z)-Ai(z)-connection}
\Bi(z)=e^{-\frac{\pi i}{6}}\Ai(ze^{-\frac{2\pi i}{3}})+e^{\frac{\pi i}{6}}\Ai(ze^{\frac{2\pi i}{3}}),
\end{equation}
and the results are
\begin{equation}\label{eq-behavior-Bi}
\left\{\begin{aligned}
\Bi(z)\sim&\frac{1}{\sqrt{\pi}}z^{-\frac{1}{4}}e^{\frac{2}{3}z^{\frac{3}{2}}}+\frac{i}{2\sqrt{\pi}}z^{-\frac{1}{4}}e^{-\frac{2}{3}z^{\frac{3}{2}}},&\arg z \in \left(-\frac{\pi}{3},\pi\right);\\
\Bi(z)\sim&\frac{1}{\sqrt{\pi}}z^{-\frac{1}{4}}e^{\frac{2}{3}z^{\frac{3}{2}}}-\frac{i}{2\sqrt{\pi}}z^{-\frac{1}{4}}e^{-\frac{2}{3}z^{\frac{3}{2}}},&\arg z \in \left(-\pi,\frac{\pi}{3}\right);\\
\Bi(z)\sim&\frac{1}{2\sqrt{\pi}}z^{-\frac{1}{4}}e^{\frac{2}{3}z^{\frac{3}{2}}}+\frac{i}{2\sqrt{\pi}}z^{-\frac{1}{4}}e^{-\frac{2}{3}z^{\frac{3}{2}}},&\arg z\in\left(\frac{\pi}{3},\frac{5\pi}{3}\right)
\end{aligned}\right.
\end{equation}
as $|z|\rightarrow\infty$. Substituting (\ref{zeta-define}) and (\ref{omega-define}) respectively into (\ref{eq-behavior-Ai}) and (\ref{eq-behavior-Bi}), we find that for $\eta$ on the Stokes line $\ell_{12}$,
\begin{equation}\label{Airy-zeta-short stokes}
\left\{\begin{aligned}
\zeta^{\frac{1}{4}}\Ai(\xi^{\frac{2}{3}}\zeta)\sim&\frac{1}{2\sqrt{\pi}}\xi^{-\frac{1}{6}}\exp\left\{-\xi\int_{\eta_{2}}^{\eta}F(\eta,\xi)^{\frac{1}{2}}ds-\xi\int_{\eta_{1}}^{\eta_{2}}F(\eta,\xi)^{\frac{1}{2}}ds\right\}, \\ \zeta^{\frac{1}{4}}\Bi(\xi^{\frac{2}{3}}\zeta)\sim&\frac{1}{\sqrt{\pi}}\xi^{-\frac{1}{6}}\exp\left\{\xi\int_{\eta_{2}}^{\eta}F(\eta,\xi)^{\frac{1}{2}}ds+\xi\int_{\eta_{1}}^{\eta_{2}}F(\eta,\xi)^{\frac{1}{2}}ds\right\}\\
&-\frac{i}{2\sqrt{\pi}}\xi^{-\frac{1}{6}}\exp\left\{-\xi\int_{\eta_{2}}^{\eta}F(\eta,\xi)^{\frac{1}{2}}ds-\xi\int_{\eta_{1}}^{\eta_{2}}F(\eta,\xi)^{\frac{1}{2}}ds\right\}
\end{aligned}\right.
\end{equation}
and
\begin{equation}\label{Airy-omega-short stokes}
\left\{\begin{aligned}
\omega^{\frac{1}{4}}\Ai(\xi^{\frac{2}{3}}\omega)&\sim\frac{1}{2\sqrt{\pi}}\xi^{-\frac{1}{6}}\exp\left\{-\xi\int_{\eta_{2}}^{\eta}F(\eta,\xi)^{\frac{1}{2}}ds\right\}, \\ \omega^{\frac{1}{4}}\Bi(\xi^{\frac{2}{3}}\omega)&\sim\frac{1}{\sqrt{\pi}}\xi^{-\frac{1}{6}}\exp\left\{\xi\int_{\eta_{2}}^{\eta}F(\eta,\xi)^{\frac{1}{2}}ds\right\}+\frac{i}{2\sqrt{\pi}}\xi^{-\frac{1}{6}}\exp\left\{-\xi\int_{\eta_{2}}^{\eta}F(\eta,\xi)^{\frac{1}{2}}ds\right\}
\end{aligned}\right.
\end{equation}
as $\xi\rightarrow+\infty$. Comparing (\ref{Airy-zeta-short stokes}) with (\ref{Airy-omega-short stokes}), we have
\begin{equation}\label{eq-zeta-omega-relation}
\zeta^{\frac{1}{4}}\left(\Ai(\xi^{\frac{2}{3}}\zeta),\Bi(\xi^{\frac{2}{3}}\zeta)\right)\sim\omega^{\frac{1}{4}}\left(\Ai(\xi^{\frac{2}{3}}\omega),\Bi(\xi^{\frac{2}{3}}\omega)\right)
\left(\begin{matrix}
e^{Q(\xi)}&-i\left[e^{Q(\xi)}+e^{-Q(\xi)}\right]\\[0.2cm]
0&e^{-Q(\xi)}
\end{matrix}\right)
\end{equation}
as $\xi\rightarrow+\infty$, uniformly for $\eta\in \gamma_{0}\cup\gamma_{1}\cup \ell_{12}$, where
\begin{equation}\label{Q(xi)}
Q(\xi)=\xi\int_{\eta_{2}}^{\eta_{1}}F(s,\xi)^{\frac{1}{2}}ds.
\end{equation}
A straightforward calculation yields (see Appendix B)
\begin{equation}\label{Q(xi)-Q0-Q1}
 Q(\xi)=i(\xi Q_{0}-Q_{1})+o(1)\quad \text{as}\quad \xi\rightarrow+\infty,
\end{equation}
where
\begin{equation}\label{eq-Q0-Q1-E0-F0-E1-F1}
Q_{0}=\frac{\sqrt{\pi}\Gamma(\frac{1}{3})}{\sqrt{3}\Gamma(\frac{11}{6})}=i(E_{0}-F_{0}), \quad Q_{1}=\frac{\pi}{3}=i(E_{1}-F_{1}).
\end{equation}

Now we begin to calculate $s_{0}$ via $\Phi_{1}(\lambda)=\Phi_{0}(\lambda)S_{0}$ and $s_{1}$ via $\Phi_{2}(\lambda)=\Phi_{1}(\lambda)S_{1}$.

If $|\eta|\rightarrow +\infty$ with $\arg\eta\sim\frac{\pi}{5}$, then $\arg\zeta\sim\frac{\pi}{3}$. Hence, substituting (\ref{zeta-eta-infty-relation}) into (\ref{eq-behavior-Ai}) and (\ref{eq-behavior-Bi}), and noting that $\lambda=\xi^{\frac{2}{5}}\eta$ and the definition of $F(\eta,\xi)$ in (\ref{eq-F(eta,xi)}), we get
\begin{equation}\label{eq-Phi-pi/5}
\begin{split}
\sqrt{2(\lambda-a)}\left(\frac{\zeta}{F(\eta,\xi)}\right)^{\frac{1}{4}}\Ai(\xi^{\frac{2}{3}}\zeta)&\sim c_{1}\frac{-1}{\sqrt{2}}\lambda^{-\frac{1}{4}}e^{-\frac{4}{5}\lambda^{\frac{5}{2}}},\\
\sqrt{2(\lambda-a)}\left(\frac{\zeta}{F(\eta,\xi)}\right)^{\frac{1}{4}}\Bi(\xi^{\frac{2}{3}}\zeta)&\sim ic_{1}\frac{-1}{\sqrt{2}}\lambda^{-\frac{1}{4}}e^{-\frac{4}{5}\lambda^{\frac{5}{2}}}+2c_{2}\frac{1}{\sqrt{2}}\lambda^{-\frac{1}{4}}e^{\frac{4}{5}\lambda^{\frac{5}{2}}}
\end{split}
\end{equation}
as $\xi, |\lambda| \rightarrow+\infty$ with $\arg\lambda\sim\frac{\pi}{5}$, where
\begin{equation}\label{eq-e1-e2}
c_{1}=-\frac{1}{\sqrt{2\pi}}\xi^{\frac{2}{15}}e^{-\xi E_{0}+E_{1}+o(1)}\quad \text{and}\quad c_{2}=\frac{1}{\sqrt{2\pi}}\xi^{\frac{2}{15}}e^{\xi E_{0}-E_{1}+o(1)}
\end{equation}
as $\xi\rightarrow+\infty$. In addition, a straightforward calculation from (\ref{eq-canonical-solutions}) leads to
\begin{equation}\label{eq-uniform-asymptotic-phik}
(\Phi_{k})_{21}\sim\frac{1}{\sqrt{2}}\lambda^{-\frac{1}{4}}e^{\frac{4}{5}\lambda^{\frac{5}{2}}}\quad\text{and}\quad (\Phi_{k})_{22}\sim\frac{-1}{\sqrt{2}}\lambda^{-\frac{1}{4}}e^{-\frac{4}{5}\lambda^{\frac{5}{2}}},\quad k\in\mathbb{Z},\quad \lambda\rightarrow\infty.
\end{equation}
Comparing (\ref{eq-Phi-pi/5}) with (\ref{eq-uniform-asymptotic-phik}), and using the results of Lemma \ref{lem-uniform-airy-large-slope}, we obtain
\begin{equation}\label{eq-Phi-pi/5-1}
\left((\Phi_{1})_{21},(\Phi_{1})_{22}\right)=\sqrt{2(\lambda-a)}\left(\frac{\zeta}{F(\eta,\xi)}\right)^{\frac{1}{4}}(\Ai(\xi^{\frac{2}{3}}\zeta),\Bi(\xi^{\frac{2}{3}}\zeta))
\left(\begin{matrix}
-\frac{i}{2c_{2}}&\frac{1}{c_{1}}\\[0.2cm]
\frac{1}{2c_{2}}&0
\end{matrix}\right)
\end{equation}
as $\xi\rightarrow+\infty$. Here, the $c_{j}$'s in (\ref{eq-Phi-pi/5-1}) are not equal but asymptotically equal to the corresponding ones in (\ref{eq-e1-e2}) as $\xi\rightarrow+\infty$. By abuse of notations, we use the same symbol for the $c_{j}$'s in these two formulas, since we  only care about the asymptotic behavior of the Stokes multipliers.

If $|\lambda|\rightarrow +\infty$ with $\lambda$ on the Stokes line $\gamma_{0}$, then $\arg\omega\sim-\frac{\pi}{3}$. Using a similar argument as in the derivation of (\ref{eq-Phi-pi/5-1}), we get
\begin{equation}\label{eq-Phi--pi/5-1}
((\Phi_{0})_{21},(\Phi_{0})_{22})=\sqrt{2(\lambda-a)}\left(\frac{\omega}{F(\eta,\xi)}\right)^{\frac{1}{4}}(\Ai(\xi^{\frac{2}{3}}\omega),\Bi(\xi^{\frac{2}{3}}\omega))
\left(\begin{matrix}
\frac{i}{2\tilde{c}_{2}}&\frac{1}{\tilde{c}_{1}}\\[0.2cm]
\frac{1}{2\tilde{c}_{2}}&0
\end{matrix}\right),
\end{equation}
where
$$\tilde{c}_{1}=-\frac{1}{\sqrt{2\pi}}\xi^{\frac{2}{15}}e^{-\xi F_{0}+F_{1}+o(1)}\quad \text{and}\quad \tilde{c}_{2}=\frac{1}{\sqrt{2\pi}}\xi^{\frac{2}{15}}e^{\xi F_{0}-F_{1}+o(1)}$$
as $\xi\rightarrow+\infty$.  Combining (\ref{eq-Phi-pi/5-1}), (\ref{eq-Phi--pi/5-1}) and (\ref{eq-zeta-omega-relation}), it is readily to get
\begin{equation}\label{eq-S0-case-I}
\begin{split}
S_{0}=&\left(\begin{matrix}\frac{i}{2\tilde{c}_{2}}&\frac{1}{\tilde{c}_{1}}\\[0.2cm]
\frac{1}{2\tilde{c}_{2}}&0\end{matrix}\right)^{-1}
\left(\begin{matrix}e^{Q(\xi)}&-i\left[e^{Q(\xi)}+e^{-Q(\xi)}\right]\\[0.2cm]
0&e^{-Q(\xi)}\end{matrix}\right)
\left(\begin{matrix}-\frac{i}{2c_{2}}&\frac{1}{c_{1}}\\[0.2cm]
\frac{1}{2c_{2}}&0\end{matrix}\right) \\[0.2cm]
=&\left(\begin{matrix}1&0\\[0.2cm]
-2i\frac{\tilde{c}_{1}}{c_{2}}\cos\{-i Q(\xi)\}&1\end{matrix}\right)
\end{split}
\end{equation}
as $\xi\rightarrow +\infty$. The last equality in (\ref{eq-S0-case-I}) holds because $Q(\xi)=\xi(F_{0}-E_{0})-(F_{1}-E_{1})+o(1)$ as $\xi\rightarrow +\infty$; cf.  (\ref{Q(xi)-Q0-Q1}) and (\ref{eq-Q0-Q1-E0-F0-E1-F1}). Finally, noting that $E_{1}+F_{1}=0$ and (\ref{Q(xi)-Q0-Q1}), we get
\begin{equation*}
\begin{split}
 s_{0}&=2ie^{-\xi (E_{0}+F_{0})+(E_{1}+F_{1})+o(1)}\cos\left\{-i Q(\xi)\right\}\\
     &=2ie^{-\xi (E_{0}+F_{0})+o(1)}\cos\left\{\xi Q_{0}-Q_{1}+o(1)\right\}
\end{split}
\end{equation*}
as $\xi\rightarrow+\infty$.

For the calculation of $s_{1}$, we need the asymptotic behaviors of the Airy functions as $|\eta|,  \xi \rightarrow +\infty$ with $\arg\eta\sim\frac{\pi}{5}\text{ and }\frac{3\pi}{5}$, {\it i.e.} $\arg\zeta\sim\frac{\pi}{3}\text{ and }\pi$, which are given in (\ref{eq-behavior-Ai}) and (\ref{eq-behavior-Bi}). Hence, following the way of deriving (\ref{eq-Phi-pi/5-1}) or (\ref{eq-Phi--pi/5-1}), we obtain
\begin{equation}\label{eq-Phi-3pi/5-1}
\begin{aligned}
((\Phi_{2})_{21},(\Phi_{2})_{22})=\sqrt{2(\lambda-a)}\left(\frac{\zeta}{F(\eta,\xi)}\right)^{\frac{1}{4}}(\Ai(\xi^{\frac{2}{3}}\zeta),\Bi(\xi^{\frac{2}{3}}\zeta))
\left(\begin{matrix}                                                                                                                                                                                          -\frac{i}{2c_{2}}&\frac{1}{2c_{1}}\\[0.2cm]                                                                                                                                                                                 \frac{1}{2c_{2}}&-\frac{i}{2c_{1}}                                                                                                                                                                       \end{matrix}\right).
\end{aligned}
\end{equation}
Combining (\ref{eq-Phi-pi/5-1}) and (\ref{eq-Phi-3pi/5-1}), one can immediately get
\begin{equation}\label{stokes-matrix-s1}
((\Phi_{2})_{21},(\Phi_{2})_{22})=((\Phi_{1})_{21},(\Phi_{1})_{22})\left(\begin{matrix}
                                                             1 & -i\frac{c_{2}}{c_{1}} \\
                                                             0 & 1
                                                           \end{matrix}\right).
\end{equation}
Noting that the constants $c_{1}$ and $c_{2}$ are defined in (\ref{eq-e1-e2}), we then have
$$s_{1}=i\exp\left\{2\xi E_{0}-2E_{1}+o(1)\right\} \quad \text{as} \quad \xi\rightarrow+\infty.$$

Finally, other Stokes multipliers can be calculated via the constraint condition (\ref{eq-constraints-stokes-multipliers}), and particularly,
\begin{equation}\label{eq-s3-large-slope}
  s_{3}=s_{-2}=i(1+s_{0}s_{1})=-i\exp\left\{2\xi \left(E_{0}-F_{0}\right)-2\left(E_{1}-F_{1}\right)+o(1)\right\}
\end{equation}
as $\xi\rightarrow+\infty$. This completes the proof of Lemma \ref{Thm-s0-large-slope}.
\qed

\subsection{Case (II): fixed $b$ and large negative $a$}
In this case, it is appropriate to make the scaling $a=-\xi^{\frac{2}{5}}$ and $\lambda=\xi^{\frac{2}{5}}\eta$ with $\xi\rightarrow+\infty$, following which, equation (\ref{Schrodinger-equation}) is reduced to
\begin{equation}\label{eq-schrodinger-y(0)-large}
\begin{split}
\frac{d^{2}Y}{d\eta^{2}}&=\xi^{2}\left[4(\eta+1)(\eta-e^{\frac{\pi i}{3}})(\eta-e^{\frac{-\pi i}{3}})+\frac{b^{2}}{\xi^{\frac{6}{5}}}+\frac{b}{\xi^{\frac{8}{5}}(\eta+1)}+\frac{3}{4\xi^{2}(\eta+1)^{2}}\right]Y\\
&:=\xi^{2}\tilde{F}(\eta,\xi)Y,
\end{split}
\end{equation}
where
\begin{equation}\label{eq-tildeF(eta,xi)}
  \tilde{F}(\eta,\xi)=4(\eta+1)(\eta-e^{\frac{\pi i}{3}})(\eta-e^{\frac{-\pi i}{3}})+\tilde{g}(\eta,\xi)
\end{equation}
and $\tilde{g}(\eta,\xi)=\mathcal{O}\left(\xi^{-\frac{6}{5}}\right)$ as $\xi\rightarrow+\infty$, uniformly for $\eta$ bounded away from $\eta=-1$. Again there are three simple turning points, say $\tilde{\eta}_{j}$, $j=0, 1, 2$ in the neighborhood of $\eta=-1$, $\eta=\alpha=e^{\frac{\pi i}{3}}$ and $\eta=\beta=e^{\frac{-\pi i}{3}}$ respectively. Moreover, we find that near the turning points $\tilde{\eta}_{1}$ and $\tilde{\eta}_{2}$, equation (\ref{eq-schrodinger-y(0)-large}) is similar to and even simpler than (\ref{eq-shrodinger-y10-large}) in Case (I), and  the Stokes geometry of $\tilde{F}(\eta,\xi)d\eta^{2}$ is the same as the one shown in Figure \ref{Figure-stokes complex-large slope}.
Hence, following the analysis in Subsection 3.1, we get
\begin{equation}
\left\{\begin{aligned}
s_{0}&=2ie^{-\xi (E_{0}+F_{0})+o(1)}\cos\left\{i R(\xi)\right\},\\
s_{1}&=ie^{2\xi E_{0}+o(1)}
\end{aligned}\right.
\end{equation}
as $\xi\rightarrow+\infty$, where $R(\xi)=\xi\int_{\tilde{\eta}_{2}}^{\tilde{\eta}_{1}}\tilde{F}(s,\xi)^{\frac{1}{2}}ds$ and $E_{0}, F_{0}$ are the same constants as the ones stated in \eqref{eq-constants-lemma-1}. Similar argument as the one of deriving (\ref{Q(xi)-Q0-Q1}) yields $R(\xi)=i \xi Q_{0}+o(1)$. Therefore
$$s_{0}=2ie^{-(-a)^{5/2}(E_{0}+F_{0})+o(1)}\cos{\left\{(-a)^{5/2}Q_{0}+o(1)\right\}}$$
as $\xi\rightarrow+\infty$.

Finally, other Stokes multipliers can be calculated via the constraint condition (\ref{eq-constraints-stokes-multipliers}), and particularly,
\begin{equation}\label{eq-s3-large-slope}
  s_{3}=s_{-2}=i(1+s_{0}s_{1})=-ie^{2\xi \left(E_{0}-F_{0}\right)+o(1)} \quad \text{as}\quad \xi\rightarrow+\infty.
\end{equation}
This completes the proof of Lemma \ref{Thm-s0-large-value}.

\subsection{Case (III): fixed $b$ and large positive $a$}
With the scaling $\lambda=\xi^{\frac{2}{5}}\eta$ and $a=\xi^{\frac{2}{5}}$, equation (\ref{Schrodinger-equation}) is reduced to
\begin{equation}\label{eq-schrodinger-large-positive-a}
\begin{split}
  \frac{d^{2}Y}{d\eta^{2}}&=\xi^{2}\left[4(\eta-1)(\eta-e^{\frac{2\pi i}{3}})(\eta-e^{\frac{-2\pi i}{3}})+\frac{b^{2}}{\xi^{\frac{6}{5}}}+\frac{b}{\xi^{\frac{8}{5}}(\eta-1)}+\frac{3}{4\xi^{2}(\eta-1)^{2}}\right]Y\\
&:=\xi^{2}\hat{F}(\eta,\xi)Y,
\end{split}
\end{equation}
where
\begin{equation}\label{eq-hatF(eta,xi)}
\hat{F}(\eta,\xi)=4(\eta-1)(\eta-e^{\frac{2\pi i}{3}})(\eta-e^{\frac{-2\pi i}{3}})+\frac{3}{4\xi^{2}(\eta-1)^{2}}+\hat{g}(\eta,\xi),\\
\end{equation}
and $(\eta-1)\hat{g}(\eta,\xi)=\mathcal{O}\left(\xi^{-\frac{6}{5}}\right)$ as $\xi\rightarrow+\infty$, uniformly for $\eta$ bounded away from $\eta=1$. Hence for large $\xi$ there are $\eta$-zeros of $\hat{F}(\eta,\xi)$, known as turning points, near $1$, $\hat{\alpha}=e^{\frac{2\pi i}{3}}$ and  $\hat{\beta}=e^{\frac{-2\pi i}{3}}$, respectively.  We note that the two turning points $\hat{\eta}_{\hat{\alpha}}$ and $\hat{\eta}_{\hat{\beta}}$ near $\hat{\alpha}$ and $\hat{\beta}$ respectively are both simple and
\begin{equation}
|\hat{\eta}_{\hat{\alpha}}-\hat{\alpha}|=\mathcal{O}\left(\xi^{-\frac{6}{5}}\right)\quad \text{and}\quad |\hat{\eta}_{\hat{\beta}}-\hat{\beta}|=\mathcal{O}\left(\xi^{-\frac{6}{5}}\right).
\end{equation}
Hence, in a neighborhood of each of these two turning points and on the stokes curves emanating form them, the Airy functions can also be used to uniformly approximate the solutions of (\ref{eq-schrodinger-large-positive-a}). Similar to (\ref{zeta-define}) and (\ref{omega-define}), if we define
\begin{equation}\label{chi-define}
  \int_{0}^{\chi}s^{\frac{1}{2}}ds=\int_{\hat{\eta}_{\hat{\alpha}}}^{\eta}\hat{F}(s,\xi)^{\frac{1}{2}}ds,
\end{equation}
then we have the following lemma which is similar to Lemma \ref{lem-uniform-airy-large-slope}.
\begin{lemma}\label{lem-uniform-airy-large-value-positive}
There are constants $D_{1}$ and $D_{2}$, depending on $\xi$, such that
\begin{equation}
Y=\left(\frac{\chi}{\hat{F}(\eta,\xi)}\right)^{\frac{1}{4}}\left\{D_{1}\left[1+o(1)\right]\Ai(\xi^{\frac{2}{3}}\chi)+D_{2}\left[1+o(1)\right]\Bi(\xi^{\frac{2}{3}}\chi)\right\},\quad \xi\rightarrow+\infty,
\end{equation}
uniformly for $\eta$ on two adjacent Stokes curves emanating from $\hat{\eta}_{\hat{\alpha}}$.
\end{lemma}
Moreover, the asymptotic behavior of $\chi(\eta)$ as $|\eta|,\, \xi\rightarrow+\infty$ is stated as follows.
\begin{lemma}\label{lemma-relation-chi-eta}
  As $\xi,|\eta|\rightarrow+\infty$ with $|\eta|\gg \xi^{2}$,
  \begin{equation}\label{eq-relation-chi-eta}
    \frac{2}{3}\chi^{\frac{3}{2}}\sim\frac{4}{5}\eta^{\frac{5}{2}}+G_{0}+o(\xi^{-1}),\qquad \arg\eta\in\left(-\frac{\pi}{5},\pi\right),
  \end{equation}
where $G_{0}=-\frac{1}{5}\B\left(\frac{1}{2},\frac{1}{6}\right)+\frac{\sqrt{3}i}{5}\B\left(\frac{1}{2},\frac{1}{6}\right)$.
\end{lemma}
The proof of this lemma is similar to that of Lemma \ref{lemma-zeta-omega-eta-infty-relation}, and hence omitted here.

Now we consider the uniform asymptotics of the solutions of (\ref{eq-schrodinger-large-positive-a}) near $\eta=1$. By a careful analysis, we find that there are three turning points $\hat{\eta}_{j}$, $j=1,2,3$, coalescing with the double pole $\eta=1$ of $\hat{F}(\eta,\xi)$ as $\xi\rightarrow+\infty$. Moreover,
\begin{equation}\label{eq-hat-eta-1-esitimate}
 \hat{\eta}_{j}-1\sim e^{\frac{j\pi i}{3}}4^{-\frac{2}{3}}\xi^{-\frac{2}{3}}\quad \text{as}\quad \xi\rightarrow+\infty.
\end{equation}
According to the idea of \cite{Dunster-1990}, near $\eta=1$, equation (\ref{eq-schrodinger-large-positive-a}) can be approximated by the following differential equation
\begin{equation}\label{eq-canonical-equation}
  \frac{d^{2}\varphi}{d\delta^{2}}=\xi^{2}\left[12\delta+\frac{3}{4\xi^{2}\delta^{2}}\right]\varphi,
\end{equation}
whose solutions can be represented by the modified Bessel functions as follows
\begin{equation}\label{eq-solution-varphi}
  \varphi_{+}(\delta)=\delta^{\frac{1}{2}}I_{\frac{2}{3}}\left(\frac{4\sqrt{3}}{3}\xi\delta^{\frac{3}{2}}\right) \quad \text{and}\quad  \varphi_{-}(\delta)=\delta^{\frac{1}{2}}K_{\frac{2}{3}}\left(\frac{4\sqrt{3}}{3}\xi\delta^{\frac{3}{2}}\right).
\end{equation}
Equation (\ref{eq-canonical-equation}) also has three turning points, say $\delta_{j}$, $j=1,2,3$, coalescing with the double pole $\delta=0$. More precisely,   $\delta_{j}\sim e^{\frac{j\pi i}{3}}4^{-\frac{2}{3}}\xi^{-\frac{2}{3}}$ for $j=1,2,3$ as $\xi\rightarrow+\infty$. Hence, from (\ref{eq-hat-eta-1-esitimate}), we see that $|\delta_{j}+1-\hat{\eta}_{j}|=o(\xi^{-\frac{2}{3}})$, $j=1,2,3$ as $\xi\rightarrow+\infty$. Define
\begin{equation}\label{eq-map-delta}
\int_{\delta_{1}}^{\delta}\left[12s+\frac{3}{4\xi^{2}s^{2}}\right]^{\frac{1}{2}}ds=\int_{\hat{\eta}_{1}}^{\eta}\hat{F}(s,\xi)^{\frac{1}{2}}ds,
\end{equation}
then $\delta(\eta)$ is a conformal mapping in the neighborhood of $\eta=1$. Moreover, we have the following lemma.

\begin{lemma}\label{lemma-uniform-bessel-large-value-positive}
There are two constants $A_{1}$ and $A_{2}$, depending on $\xi$, such that
\begin{equation}
Y=\left[\frac{12\delta+\frac{3}{4\xi^{2}\delta^{2}}}{\hat{F}(\eta,\xi)}\right]^{\frac{1}{4}}\left\{A_{1}\left[1+o(1)\right]\varphi_{+}(\delta)+A_{2}\left[1+o(1)\right]\varphi_{-}(\delta)\right\}
\end{equation}
uniformly for $\eta$ on any two adjacent Stokes lines emanating from $\hat{\eta}_{j}$, $j=1,2,3$.
\end{lemma}
The proof of this lemma is similar to those of \cite[Theorems 1 and 2]{APC}.
By evaluating the integrals in (\ref{eq-map-delta}) (see Appendix C), we get the asymptotic behavior of $\delta(\eta)$ as follows.
\begin{lemma}\label{lemma-relation-eta-delta}
As $\xi, |\eta| \rightarrow+\infty$ with $|\eta|\gg\xi^{2}$, the asymptotic behavior of $\delta(\eta)$ is given by
  \begin{equation}\label{eq-relation-eta-delta}
    \frac{4\sqrt{3}}{3}\delta^{\frac{3}{2}}=\frac{4}{5}\eta^\frac{5}{2}-H_{0}+o(\xi^{-1}),\qquad \arg\eta\in\left(-\frac{3\pi}{5},\frac{3\pi}{5}\right),
  \end{equation}
where $H_{0}=\frac{2}{5}\B\left(\frac{1}{6},\frac{1}{2}\right)$.
\end{lemma}

Now we turn to the proof of Lemma \ref{Thm-large-initial-value-postive}, {\it i.e.}, to derive the asymptotic behavior of the Stokes multipliers $s_{k}$ as $\xi\rightarrow+\infty$, $k=0,1,2,3,4$.

\textbf{Proof of Lemma \ref{Thm-large-initial-value-postive}.} To calculate $s_{0}$, we use the uniform asymptotics of the solutions of (\ref{eq-schrodinger-large-positive-a}) near $\eta=1$; see Lemmas \ref{lemma-uniform-bessel-large-value-positive} and \ref{lemma-relation-eta-delta}. If $\arg\eta\sim\frac{\pi}{5}$ as $|\eta|\rightarrow\infty$, then $\arg\delta^{\frac{3}{2}}\sim\frac{\pi}{2}$.
Recall the uniform asymptotic behaviors of the modified Bessel functions stated in \cite[(10.4.2), (10.40.5)]{NIST-handbook} and $\lambda=\xi^{\frac{2}{5}}\eta$. Then (\ref{eq-solution-varphi}) and (\ref{eq-relation-eta-delta}) lead to
\begin{equation}\label{eq-modified-bessel-function-pi/2}
\left\{\begin{aligned}
\sqrt{2(\lambda-a)}\left[\frac{12\delta+\frac{3}{4\xi^{2}\delta^{2}}}{\hat{F}(\eta,\xi)}\right]^{\frac{1}{4}}\varphi_{+}(\delta)\sim& f_{1}\frac{1}{\sqrt{2}\lambda^{\frac{1}{4}}}e^{\frac{4}{5}\lambda^{\frac{5}{2}}}+f_{2}\frac{-1}{\sqrt{2}\lambda^{\frac{1}{4}}}e^{-\frac{4}{5}\lambda^{\frac{5}{2}}},\\
\sqrt{2(\lambda-a)}\left[\frac{12\delta+\frac{3}{4\xi^{2}\delta^{2}}}{\hat{F}(\eta,\xi)}\right]^{\frac{1}{4}}\varphi_{-}(\delta)\sim& f_{3}\frac{-1}{\sqrt{2}\lambda^{\frac{1}{4}}}e^{-\frac{4}{5}\lambda^{\frac{5}{2}}}
\end{aligned}\right.
\end{equation}
as $|\eta|,\xi\rightarrow+\infty$ with $\arg\eta\sim\frac{\pi}{5}$, where
\begin{equation}\label{eq-f1-f2-f3}
f_{1}=\frac{3^{1/4}}{2\sqrt{\pi}}\xi^{-\frac{1}{5}}e^{-\xi H_{0}+o(1)}\quad \text{and}\quad f_{2}=-\frac{i3^{1/4}e^{\frac{2\pi i}{3}}}{2\sqrt{\pi}}\xi^{-\frac{1}{5}}e^{\xi H_{0}+o(1)}=\frac{e^{\frac{2\pi i}{3}}}{\pi}f_{3},
\end{equation}
as $\xi\rightarrow+\infty$. Comparing (\ref{eq-modified-bessel-function-pi/2}) with (\ref{eq-uniform-asymptotic-phik}), we get
\begin{equation}\label{eq-Sk-bessel-function-pi/5}
\left((\Phi_{1})_{21},(\Phi_{1})_{22}\right)=\sqrt{2(\lambda-a)}\left[\frac{4\delta+\frac{3}{4\xi^{2}\delta^{2}}}{\hat{F}(\eta,\xi)}\right]^{\frac{1}{4}}\left(\varphi_{+}(\delta),\varphi_{-}(\delta)\right)
\left(\begin{matrix}
  \frac{1}{f_{1}} & 0 \\[0.2cm]
  -\frac{f_{2}}{f_{1}f_{3}} & \frac{1}{f_{3}}
\end{matrix}\right).
\end{equation}
If $\arg\eta\sim\frac{-\pi}{5}$ as $|\eta|\rightarrow +\infty$, then $\arg\delta^{\frac{3}{2}}\sim\frac{-\pi}{2}$. Using a similar argument for deriving (\ref{eq-Sk-bessel-function-pi/5}), one has
\begin{equation}\label{eq-Sk-bessel-function--pi/5}
\left((\Phi_{0})_{21},(\Phi_{0})_{22}\right)=\sqrt{2(\lambda-a)}\left[\frac{4\delta+\frac{3}{4\xi^{2}\delta^{2}}}{\hat{F}(\eta,\xi)}\right]^{\frac{1}{4}}\left(\varphi_{+}(\delta),\varphi_{-}(\delta)\right)
\left(\begin{matrix}
  \frac{1}{\hat{f}_{1}} & 0 \\[0.2cm]
  -\frac{\hat{f}_{2}}{\hat{f}_{1}\hat{f}_{3}} & \frac{1}{\hat{f}_{3}}
\end{matrix}\right),
\end{equation}
where $\hat{f}_{1}=f_{1}$, $\hat{f}_{2}=-e^{\frac{-4\pi i}{3}}f_{2}$ and $\hat{f}_{3}=f_{3}$.
Combining (\ref{eq-Sk-bessel-function-pi/5}) and (\ref{eq-Sk-bessel-function--pi/5}), we have
\begin{equation}\label{eq-S0-large-positive-a}
  S_{0}=\left(\begin{matrix}
                1 & 0 \\
                \frac{\hat{f}_{2}-f_{2}}{f_{1}} & 1
              \end{matrix}\right).
\end{equation}
Therefore $s_{0}=\frac{\hat{f}_{2}-f_{2}}{f_{1}}=-ie^{2\xi H_{0}+o(1)}$ as $\xi\rightarrow+\infty$.

Next, to derive $s_{1}$, we need the uniform asymptotics of the solutions of (\ref{eq-schrodinger-large-positive-a}) near $\eta=\hat{\alpha}$ (see Lemmas \ref{lem-uniform-airy-large-value-positive} and \ref{lemma-relation-chi-eta}). Following the analysis of deriving $s_{1}$ in case (I) or (II), we obtain
\begin{equation}\label{eq-s1-large-value-positive}
  s_{1}=i\exp\left\{2\xi G_{0}+o(1)\right\}, \quad \text{as}\quad \xi\rightarrow+\infty.
\end{equation}

Other Stokes multipliers can be derived by the constraint condition (\ref{eq-constraints-stokes-multipliers}) and (\ref{eq-sk-s-k-relation}), and particularly,
\begin{equation}\label{eq-s3-large-slope}
  s_{3}=s_{-2}=i(1+s_{0}s_{1})=i\exp\left\{2\xi \left(G_{0}+H_{0}\right)+o(1)\right\}
\end{equation}
and
\begin{equation}\label{eq-s3-large-slope}
  s_{2}=-\overline{s_{3}}=i\exp\left\{-2\xi G_{0}+o(1)\right\}
\end{equation}
as $\xi\rightarrow+\infty$. This completes the proof of Lemma \ref{Thm-large-initial-value-postive}.
\qed

\section{Discussion}\label{sec:discussion}
We have considered the initial value problem of (\ref{PI equation}) with initial data (\ref{PI-equation-initial-problem}) in three special cases. In fact we have given a rigorous proof of the conclusions in \cite{Bender-Komijani-2015} for PI, built more precise relations between the initial value of PI solutions and their large negative $t$ asymptotic behavior, and obtained the limiting form connection formulas (\ref{eq-limit-connection-h})--(\ref{eq-limit-connection-rho-sigma-case-II}). There are still several issues to be further investigated.


First, we have only considered three special cases of the initial value problem of PI. An equally natural question is   what   would happen when both $a$ and $b$ are large.
More attention should be paid to \cite[Fig.\,4.5]{Fornberg-Weideman} by  Fornberg and Weideman.
According to this figure, we find that it may be divided into two cases:
\begin{equation}\label{eq-large-a-b}
\begin{array}{ll}
  \text{(1)} &  \text{large $|b|$ and large negative $a$},\\
  \text{(2)} &  \text{large $|b|$ and large positive $a$}.
\end{array}
\end{equation}
Moreover, we see that the analysis of the first case in (\ref{eq-large-a-b}) is similar to  one of case (I) or case (II) in the present paper, while for the second case in (\ref{eq-large-a-b}), more careful analysis is needed. For example, a description of the Stokes geometry and the Stokes curves at the finite plane seems to be vital; cf. \cite{Kawai-Takei-2005-book}.



Of course, one may consider other special cases of Clarkson's open problem. For example $y(0)=y'(0)=0$. According to \cite{Joshi-Kruskal-1992}, see also in \cite{Qin-Lu-2008}, the PI solution with $y(0)=y'(0)=0$  oscillates about $y=-\sqrt{\frac{-t}{6}}$ when $t<0$,  and it satisfies (\ref{oscillate-behavior-t-infty}) and (\ref{phi(t)}). Hence, a question is how to determine the exact values of the parameters $d$ and $\theta$ in (\ref{oscillate-behavior-t-infty}) and (\ref{phi(t)}). In fact, the result is already known, and is given by Kitaev \cite{Kitaev-1995} via the WKB method. To the best of our knowledge, it may be the only special case of the Clarkson's open problem which has been fully solved for PI equation.

However, a more important  and challenging problem is the general case of Clarkson's open problem. Assuming that  $a$ and $b$ are fixed parameters, then equation (\ref{Schrodinger-equation}) has a regular singularity at $\lambda=a$ and an irregular singularity at $\lambda=\infty$ of rank 3. As far as we know, there is no such known special function that can be used to approximate the solutions   of (\ref{Schrodinger-equation}), uniformly for $\lambda$ near both $a$ and $\infty$. This is the main difficulty of the general case of Clarkson's open problem.

Finally, it seems quite promising that  the method of {\it uniform asymptotics} can also be applied to   the initial value problems of other Painlev\'{e} equations. In fact, similar results for the properties of the  PII  solutions have also been stated in \cite{Bender-Komijani-2015}.
\vskip 1cm

\section*{Acknowledgements}
The work of Wen-Gao Long was supported in part by the National Natural Science Foundation of China under Grant Number 11571376 and GuangDong Natural Science Foundation under Grant Number 2014A030313176.
The work of Yu-Tian Li was supported in part by the Hong Kong Research Grants Council
[grant numbers 201513 and  12303515] and the HKBU Strategic Development Fund.
The work of Sai-Yu Liu was supported in part by   the National Natural Science Foundation of China under grant numbers  11326082 and 11401200.
The work of   Yu-Qiu Zhao  was supported in part by the National
Natural Science Foundation of China under grant numbers
10871212 and  11571375.

\vskip 1cm

\begin{appendix}

\setcounter{equation}{0}
\renewcommand{\theequation}{A.\arabic{equation}}

\section*{Appendix A:  Proof of Lemma \ref{lemma-zeta-omega-eta-infty-relation}}
 The idea to prove Lemma \ref{lemma-zeta-omega-eta-infty-relation} is to derive the asymptotic behavior of the integral on the right hand side of (\ref{zeta-define}) as $\xi, |\eta|\rightarrow+\infty$. Here we only show the validity of (\ref{zeta-eta-infty-relation}). The proof of (\ref{omega-eta-infty-relation}) is similar and hence omitted here.

Choose $T^{\star}$ to be a point such that $|T^{\star}-\eta_{1}|\sim\xi^{-\frac{3}{4}}$, then
\begin{equation}\label{eq-appendix-1}
\begin{split}
\int_{\eta_{1}}^{\eta}F(s,\xi)^{\frac{1}{2}}ds=\int_{\eta_{1}}^{T^{\star}}F(s,\xi)^{\frac{1}{2}}ds+\int_{T^{\star}}^{\eta}F(s,\xi)^{\frac{1}{2}}ds:=I_{1}+I_{2}.
\end{split}
\end{equation}
Noting that $F(\eta,\xi)=\mathcal{O}(\xi^{-\frac{3}{4}})$ as $\xi\rightarrow+\infty$ uniformly for $\eta$ on the integration contour from $\eta_{1}$ to $T^{\star}$, we immediately conclude that $I_{1}=o(\xi^{-1})$. Since $|T^{\star}-\eta_{1}|\sim\xi^{-\frac{3}{4}}$ and $\eta_{1}-\alpha\sim-\frac{1}{6\xi}$ as $\xi\rightarrow+\infty$ (see the paragraph under (\ref{eq-F(eta,xi)})), then $|T^{\star}-\alpha|>c \, \xi^{-\frac{3}{4}}$ for some constant $c>0$. By the Taylor expansion of $F(s,\xi)$ with respect to $\frac{1}{\xi}$, we have
\begin{equation}\label{eq-F-taylor-expansion-T*-eta}
\begin{aligned}
F(s,\xi)^{\frac{1}{2}}=&2(s+1)^{\frac{1}{2}}(s-\alpha)^{\frac{1}{2}}(s-\beta)^{\frac{1}{2}}-\frac{1}{2\xi s(s+1)^{\frac{1}{2}}(s-\alpha)^{\frac{1}{2}}(s-\beta)^{\frac{1}{2}}}\\
&+\mathcal{O}\left(\frac{1}{\xi^{\frac{6}{5}}\sqrt{(s+1)(s-\alpha)(s-\beta)}}\right)
\end{aligned}
\end{equation}
as $\xi\rightarrow+\infty$, uniformly for $s$ on the integration contour from $T^{\star}$ to $\eta$, where $\alpha=\bar{\beta}=e^{\frac{\pi i}{3}}$. Therefore,
\begin{equation*}
\begin{split}
I_{2}&=\int_{T^{\star}}^{\eta}2(s+1)^{\frac{1}{2}}(s-\alpha)^{\frac{1}{2}}(s-\beta)^{\frac{1}{2}}ds-\int_{T^{\star}}^{\eta}\frac{1}{2\xi s(s+1)^{\frac{1}{2}}(s-\alpha)^{\frac{1}{2}}(s-\beta)^{\frac{1}{2}}}ds+o\left(\frac{1}{\xi}\right)\\
&=\int_{T^{\star}}^{\eta}2(s+1)^{\frac{1}{2}}(s-\alpha)^{\frac{1}{2}}(s-\beta)^{\frac{1}{2}}ds-\int_{T^{\star}}^{\infty}\frac{1}{2\xi s(s+1)^{\frac{1}{2}}(s-\alpha)^{\frac{1}{2}}(s-\beta)^{\frac{1}{2}}}ds+o\left(\frac{1}{\xi}\right)\\
\end{split}
\end{equation*}
provided  that $|\eta|\gg 1$ as $\xi\rightarrow+\infty$. Recalling the facts that $|T^{\star}-\eta_{1}|\sim\xi^{-\frac{3}{4}}$ and $\eta_{1}-\alpha\sim -\frac{1}{6\xi}$, we further obtain
\begin{equation}\label{eq-appendix-2}
\begin{split}
I_{2}=\int_{\alpha}^{\eta}2(s+1)^{\frac{1}{2}}(s-\alpha)^{\frac{1}{2}}(s-\beta)^{\frac{1}{2}}ds-\frac{E_{1}}{\xi}+o\left(\frac{1}{\xi}\right),
\end{split}
\end{equation}
where
\begin{equation}\label{eq-E1-appendix}
E_{1}=\int_{\alpha}^{\infty}\frac{1}{2s(s+1)^{\frac{1}{2}}(s-\alpha)^{\frac{1}{2}}(s-\beta)^{\frac{1}{2}}}ds=-\frac{\pi i}{6}.
\end{equation}
Next, using integration by parts, we have
\begin{equation}\label{eq-appendix-3}
2\int_{\alpha}^{\eta}(s+1)^{\frac{1}{2}}(s-\alpha)^{\frac{1}{2}}(s-\beta)^{\frac{1}{2}}ds=E_{0}+\frac{4}{5}\eta^{\frac{5}{2}}+o\left(\frac{1}{\xi}\right)
\end{equation}
as $\xi, |\eta|\rightarrow +\infty$ provided that $|\eta|\gg \xi^{2}$ and $\arg\eta\in\left(-\frac{\pi}{5},\pi\right)$, where
\begin{equation}\label{eq-E0-appendix}
  E_{0}=\frac{6}{5}\int_{\alpha}^{\infty}\frac{1}{\sqrt{s^3+1}}ds=\frac{3}{5}\B\left(\frac{1}{2},\frac{1}{3}\right)-\frac{\sqrt{3}i}{5}\B\left(\frac{1}{2},\frac{1}{3}\right).
\end{equation}
Finally, combining (\ref{eq-appendix-1}), (\ref{eq-appendix-2}) and (\ref{eq-appendix-3}), we obtain the desired formula (\ref{zeta-eta-infty-relation}).

\setcounter{equation}{0}
\renewcommand{\theequation}{B.\arabic{equation}}

\section*{Appendix B:  Proof of (\ref{Q(xi)-Q0-Q1}) and (\ref{eq-Q0-Q1-E0-F0-E1-F1})}
First, similar to the proof of Lemma \ref{lemma-zeta-omega-eta-infty-relation}, we choose $T_{1}^{\star}$ and $T_{2}^{\star}$ such that $|T_{j}^{\star}-\eta_{j}|\sim\xi^{-\frac{3}{4}}$ for $j=1,2$  as $\xi\rightarrow+\infty$. Then we split the integral in (\ref{Q(xi)}) to give
\begin{equation}\label{eq-Q(xi)-three-integral}
  Q(\xi)=\xi\int_{\eta_{2}}^{T_{2}^{\star}}F(s,\xi)ds+\xi\int_{T_{2}^{\star}}^{T_{1}^{\star}}F(s,\xi)ds+\xi\int_{T_{1}^{\star}}^{\eta_{1}}F(s,\xi)ds:=K_{1}+K_{2}+K_{3}.
\end{equation}
By the same argument of estimating $I_{1}$ in Appendix A, we have $K_{1},K_{3}=o(1)$ as $\xi\rightarrow+\infty$. Next, we estimate $K_{2}$.
Making use of the facts $\eta_{1}-\alpha\sim-\frac{1}{6\xi}$,  $\eta_{2}-\beta\sim-\frac{1}{6\xi}$ (see the paragraph under (\ref{eq-F(eta,xi)})) and $|T_{j}^{\star}-\eta_{j}|\sim\xi^{-\frac{3}{4}}$ as $\xi\rightarrow +\infty$ , we conclude that
\begin{equation}\label{eq-Tj-alpha-beta-estimate}
|T_{1}^{\star}-\alpha|\sim\xi^{-\frac{3}{4}} \quad \text{and}\quad |T_{2}^{\star}-\beta|\sim\xi^{-\frac{3}{4}}.
\end{equation}
Then, by the Taylor expansion of $F(s,\xi)$ with respect to $\frac{1}{\xi}$, we get
\begin{equation}\label{eq-F-taylor-expansion-T*-eta}
\begin{aligned}
F(s,\xi)^{\frac{1}{2}}=&2(s+1)^{\frac{1}{2}}(s-\alpha)^{\frac{1}{2}}(s-\beta)^{\frac{1}{2}}+\frac{1}{2\xi s(s+1)^{\frac{1}{2}}(s-\alpha)^{\frac{1}{2}}(s-\beta)^{\frac{1}{2}}}\\
&+\mathcal{O}\left(\frac{1}{\xi^{\frac{6}{5}}\sqrt{(s+1)(s-\alpha)(s-\beta)}}\right)
\end{aligned}
\end{equation}
as $\xi\rightarrow+\infty$ uniformly for $s$ on the integration contour from $T_{2}^{\star}$ to $T_{1}^{\star
}$. Hence,
\begin{equation*}
\begin{split}
K_{2}=\xi\int_{T_{2}^{\star}}^{T_{1}^{\star}}F(s,\xi)^{\frac{1}{2}}ds
=&2\xi\int_{T_{2}^{\star}}^{T_{1}^{\star}}(s+1)^{\frac{1}{2}}(s-\alpha)^{\frac{1}{2}}(s-\beta)^{\frac{1}{2}}ds\\
&-\frac{1}{2}\int_{T_{2}^{\star}}^{T_{1}^{\star}}\frac{1}{(s+1)^{\frac{1}{2}}(s-\alpha)^{\frac{1}{2}}(s-\beta)^{\frac{1}{2}}s}ds+o(1).\\
\end{split}
\end{equation*}
Using (\ref{eq-Tj-alpha-beta-estimate}) again, we can further obtain
\begin{equation}\label{eq-estimation-K2}
\begin{split}
K_{2}=&2\xi\int_{\beta}^{\alpha}(s+1)^{\frac{1}{2}}(s-\alpha)^{\frac{1}{2}}(s-\beta)^{\frac{1}{2}}ds\\
&-\frac{1}{2}\int_{\beta}^{\alpha}\frac{1}{(s+1)^{\frac{1}{2}}(s-\alpha)^{\frac{1}{2}}(s-\beta)^{\frac{1}{2}}s}ds+o(1)
\end{split}
\end{equation}
as $\xi\rightarrow+\infty$. Therefore, combining (\ref{eq-Q(xi)-three-integral}) and (\ref{eq-estimation-K2}), we get
\begin{equation}\label{Q(xi)-Q0-Q1-appendix}
 Q(\xi)=i(\xi Q_{0}-Q_{1})+o(1)\quad \text{as}\quad \xi\rightarrow+\infty,
\end{equation}
where
\begin{equation}\label{eq-Q0-Q1-E0-F0-E1-F1-appendix}
\begin{split}
Q_{0}&=-2i\int_{\beta}^{\alpha}\sqrt{(s+1)(s-\alpha)(s-\beta)}ds=\frac{\sqrt{\pi}\Gamma(\frac{1}{3})}{\sqrt{3}\Gamma(\frac{11}{6})}=\frac{2\sqrt{3}}{5}\B\left(\frac{1}{2},\frac{1}{3}\right), \\
Q_{1}&=-i\int_{\beta}^{\alpha}\frac{ds}{2\sqrt{(s+1)(s-\beta)(s-\alpha)}s}=\frac{\pi}{3}.
\end{split}
\end{equation}
Finally, in view of  the explicit representations (\ref{eq-constants-lemma-1}) of $E_{0}, F_{0}, E_{1}$ and $F_{1}$, we see that
$F_{1}-E_{1}=iQ_{1}$ and $F_{0}-E_{0}=iQ_{0}$.

\setcounter{equation}{0}
\renewcommand{\theequation}{C.\arabic{equation}}

\section*{Appendix C:  Proof of Lemma \ref{lemma-relation-eta-delta}}
The idea to prove Lemma \ref{lemma-relation-eta-delta} is to derive an asymptotic approximation for each    integral  in (\ref{eq-map-delta}).

Choose $\delta^{\star}$ to be a point with $|\delta^{\star}-\delta_{1}|\sim \xi^{-\frac{1}{2}}$ as $\xi\rightarrow+\infty$, then
\begin{equation}\label{eq-integral-delta}
\int_{\delta_{1}}^{\delta}\left[12s+\frac{3}{4\xi^{2}s^{2}}\right]^{\frac{1}{2}}ds=\int_{\delta_{1}}^{\delta^{\star}}\left[12s+\frac{3}{4\xi^{2}s^{2}}\right]^{\frac{1}{2}}ds+\int_{\delta^{\star}}^{\delta}\left[12s+\frac{3}{4\xi^{2}s^{2}}\right]^{\frac{1}{2}}ds:=J_{1}+J_{2}.
\end{equation}
Similarly, we split the integral on the right-hand side of (\ref{eq-map-delta}) as follows
\begin{equation}\label{eq-integral-eta-F}
 \int_{\hat{\eta}_{1}}^{\eta}\hat{F}(s,\xi)^{\frac{1}{2}}ds=\int_{\hat{\eta}_{1}}^{1+\delta^{\star}}\hat{F}(s,\xi)^{\frac{1}{2}}ds+\int_{1+\delta^{\star}}^{\eta}\hat{F}(s,\xi)^{\frac{1}{2}}ds:=\hat{J}_{1}+\hat{J}_{2}.
\end{equation}
Since when $s$ is on the integration contour from $\eta_{1}$ to $1+\delta^{\star}$,
\begin{equation}\label{eq-doubletilde-F-quyu-eta=1}
\begin{aligned}
\left[\hat{F}(s,\xi)\right]^{\frac{1}{2}}&=\left[12(s-1)+\frac{3}{4\xi^{2}(s-1)^{2}}+\mathcal{O}(\xi^{-\frac{14}{15}})\right]^{\frac{1}{2}}\\
&=\left[12(s-1)+\frac{3}{4\xi^{2}(s-1)^{2}}\right]^{\frac{1}{2}}+\mathcal{O}(\xi^{-\frac{3}{5}})
\end{aligned}
\end{equation}
as $\xi\rightarrow+\infty$, then
\begin{equation}\label{eq-integral-tildeJ1}
  \hat{J}_{1}=\int_{\hat{\eta}_{1}-1}^{\delta^{\star}}\left[12s+\frac{3}{4\xi^{2}s^{2}}\right]^{\frac{1}{2}}ds+\mathcal{O}(\xi^{-\frac{11}{10}})=J_{1}+o(\xi^{-1})
\end{equation}
as $\xi\rightarrow+\infty$. Here, we have made use of the fact that $|\hat{\eta}_{1}-1-\delta_{1}|=o(\xi^{-\frac{2}{3}})$, mentioned in the paragraph above (\ref{eq-map-delta}). Combining (\ref{eq-integral-delta}), (\ref{eq-integral-eta-F}) and (\ref{eq-integral-tildeJ1}), we find that (\ref{eq-map-delta}) is reduced to
\begin{equation}\label{eq-J2-tildeJ2-relation}
  J_{2}=\hat{J}_{2}+o(\xi^{-1}) \qquad \text{as}\quad \xi\rightarrow+\infty.
\end{equation}

Now we turn to the estimation of $J_{2}$ and $\hat{J}_{2}$. Since $\frac{1}{\xi^{2}s^{3}}=o(1)$ as $\xi\rightarrow+\infty$, uniformly for $s$ on the integration contour from $\delta^{\star}$ to $\delta$, then
\begin{equation}\label{eq-delta-taylor-expansion}
  \left[12s+\frac{3}{4\xi^{2}s^{2}}\right]^{\frac{1}{2}}=2\sqrt{3}s^{\frac{1}{2}}+\mathcal{O}\left(\frac{1}{\xi^{2}s^{\frac{5}{2}}}\right)\quad \text{as}\quad \xi\rightarrow+\infty.
\end{equation}
Hence, we have
\begin{equation}\label{eq-J2-approximation}
 J_{2}=\frac{4\sqrt{3}}{3}\delta^{\frac{3}{2}}-\frac{4\sqrt{3}}{3}(\delta^{\star})^{\frac{3}{2}}+o(\xi^{-1})  \quad \text{as}\quad \xi\rightarrow+\infty.
\end{equation}
Similarly, when $s$ is on the integration contour from $\delta^{\star}+1$ to $\eta$,
\begin{equation}\label{eq-eta-taylor-expansion}
  \hat{F}(s,\xi)^{\frac{1}{2}}=2(s-1)^{\frac{1}{2}}(s-\tilde{\alpha})^{1/2}(s-\tilde{\beta})^{\frac{1}{2}}+\mathcal{O}\left(\frac{1}{\xi^{2}(s-1)^{\frac{5}{2}}}\right)\quad \text{as}\quad \xi\rightarrow+\infty,
\end{equation}
since $|\delta^{\star}|\sim \xi^{-\frac{1}{2}}$, as can be derived from the facts $\delta_{1}\sim e^{\frac{j\pi i}{3}}4^{-\frac{2}{3}}\xi^{-\frac{2}{3}}$ (see the paragraph under (\ref{eq-solution-varphi})) and $|\delta^{\star}-\delta_{1}|\sim\xi^{-\frac{1}{2}}$.
Then, we have
\begin{equation}\label{eq-tildeJ2-approximation}
\begin{aligned}
\hat{J}_{2}&=2\int_{\delta^{\star}+1}^{\eta}(s-1)^{\frac{1}{2}}(s-\tilde{\alpha})^{\frac{1}{2}}(s-\tilde{\beta})^{\frac{1}{2}}ds+o(\xi^{-1})\\
           &=2\int_{1}^{\eta}(s^{3}-1)^{\frac{1}{2}}ds-2\int_{1}^{\delta^{\star}+1}(s^{3}-1)^{\frac{1}{2}}ds+o\left(\xi^{-1}\right)
\end{aligned}
\end{equation}
as $\xi\rightarrow+\infty$. Using the fact $|\delta^{\star}|\sim \xi^{-\frac{1}{2}}$ again, a simple calculation yields
$$2\int_{1}^{\delta^{\star}+1}(s^{3}-1)^{\frac{1}{2}}ds=\frac{4\sqrt{3}}{3}\left(\delta^{\star}\right)^{\frac{3}{2}}+o\left(\xi^{-1}\right),\quad \text{as}\quad \xi\rightarrow\infty.$$
Therefore,
\begin{equation}\label{eq-tildeJ2-leading-behavior}
\begin{aligned}
\hat{J}_{2}=&2\int_{1}^{\eta}(s^{3}-1)^{\frac{1}{2}}ds-\frac{4\sqrt{3}}{3}\left(\delta^{\star}\right)^{\frac{3}{2}}+o\left(\xi^{-1}\right)\\
 =&\frac{4}{5}\eta^{\frac{5}{2}}-H_{0}-\frac{4\sqrt{3}}{3}\left(\delta^{\star}\right)^{\frac{3}{2}}+o\left(\xi^{-1}\right)
\end{aligned}
\end{equation}
as $\xi, |\eta|\rightarrow+\infty$, provided that $|\eta|\gg \xi^{2}$ and $\arg\eta\in\left(-\frac{3\pi}{5},\frac{3\pi}{5}\right)$, where
\begin{equation}\label{eq-G0-appendix}
 H_{0}=\frac{6}{5}\int_{1}^{\infty}\frac{1}{\sqrt{s^{3}-1}}ds=\frac{2}{5}\B\left(\frac{1}{6},\frac{1}{2}\right).
\end{equation}
A combination of (\ref{eq-J2-tildeJ2-relation}), (\ref{eq-J2-approximation}) and (\ref{eq-tildeJ2-leading-behavior}) leads to (\ref{eq-relation-eta-delta}). This completes the proof of Lemma \ref{lemma-relation-eta-delta}.

\setcounter{equation}{0}
\renewcommand{\theequation}{D.\arabic{equation}}

\section*{Appendix D:  Calculation of integrals in \eqref{eq-E0-appendix} and \eqref{eq-Q0-Q1-E0-F0-E1-F1-appendix}}

To verify  the first formula  in \eqref{eq-Q0-Q1-E0-F0-E1-F1-appendix}, it suffices to show that
\begin{equation}\label{evaluation-Q0}
I=\int^\alpha_\beta \sqrt{(s+1)(s-\alpha)(s-\beta)} ds=\frac {i\sqrt 3 } 5 B\left (\frac 1 2, \frac 1 3\right ),
\end{equation}
where $\alpha=e^{\pi i/3}$, $\beta=e^{-\pi i/3}$, the principal branches in the integrand are chosen such that $\arg(s+1),~\arg(s-\alpha),~\arg(s-\beta)\in (-\pi, \pi)$, and the integration path is  the line segment connecting $\beta$ and $\alpha$, subject to deformation.

It is readily seen that  $I= \int^\alpha_\beta \sqrt{ s^3+1} \; ds$. Integrating by parts once, we have
\begin{equation*}
 I=-\frac 3 2 I +\frac 3 2\int^\alpha_\beta \frac { ds}{  \sqrt{ s^3+1}}.
\end{equation*}
Applying Cauchy's integral theorem,  we further have
\begin{equation*}
I=\frac 3 5 \int^\alpha_\beta \frac { ds}{  \sqrt{ s^3+1}}=\frac 3 5\left [ \int^{\infty e^{-\pi i/3}}_\beta  \frac { ds}{  \sqrt{ s^3+1}}  -\int^{\infty e^{\pi i/3}}_\alpha  \frac { ds}{  \sqrt{ s^3+1}}\right ],
\end{equation*}
where, in the last integrals, we may take the integration paths to be portions of the rays emanating   from the origin. Now   a change of variable $t=se^{\pm \pi i/3}$ gives
$$I=\frac {3 i} 5 \int^{+\infty}_1 \frac {dt}{\sqrt{t^3-1}}=\frac i 5 B\left (\frac 1 2, \frac 1 6\right ) =\frac {i\sqrt 3} 5 B\left (\frac 1 2, \frac 1 3\right ).$$
Here, in the last equality, use has been made of the reflection formula of the gamma function.  Thus we have \eqref{evaluation-Q0}, and the formula in \eqref{eq-Q0-Q1-E0-F0-E1-F1-appendix} for $Q_0$ follows accordingly.

The expression  in \eqref{eq-Q0-Q1-E0-F0-E1-F1-appendix}  for $Q_1$ can be derived similarly, and probably  easier. No integration by parts  is needed in this case.

Other integrals can also be evaluated. For instance, for the quantity  $E_0$ in \eqref{eq-E0-appendix}, we have
 \begin{equation*}
 E_0=\frac 6 5 \int^{\infty e^{\pi i/3}}_\alpha \frac{ ds}{  \sqrt{ s^3+1}}= \frac 6 5 \left (-i e^{\frac {\pi i} 3}\right ) \int^{+\infty}_1 \frac {dt}{\sqrt{t^3-1}}=\frac {\sqrt 3 (\sqrt 3 -i)} 5 B\left (\frac 1 2, \frac 1 3\right ).
 \end{equation*}

\end{appendix}

\end{document}